\documentclass[a4paper, 12pt]{amsart}

\usepackage[T1]{fontenc}
\usepackage[utf8]{inputenc}
\usepackage[english]{babel}

\usepackage{tikz}
\usepackage{tikz-cd}
\usetikzlibrary[arrows.meta]
\usetikzlibrary{arrows,decorations.pathreplacing,calligraphy,patterns,positioning,matrix,calc,bending,decorations.pathmorphing,shapes.geometric}
\tikzcdset{
  scale cd/.style={
    every label/.append style={scale=#1},
    cells={nodes={scale=#1}}
  }}

\usepackage{amsfonts}
\usepackage{amssymb}
\usepackage{amsmath}
\usepackage{amsthm}

\usepackage{enumitem}
\usepackage{graphicx}
\usepackage{subcaption}
\usepackage{moresize} 

\theoremstyle{definition}
\newtheorem{de}{Definition}[section]
\theoremstyle{plain}
\newtheorem{thm}[de]{Theorem}
\newtheorem{mainthm}{Theorem}

\newtheorem{lem}[de]{Lemma}
\newtheorem{p}[de]{Proposition}
\newtheorem{cor}[de]{Corollary}

\theoremstyle{definition}

\usepackage{newcommands}

\RequirePackage{doi}
\usepackage{hyperref}
\usepackage[libertine]{newtxmath}
\usepackage[ttscale=.875]{libertine}

\title[Determinant of the distance matrix of a tree]{Non-intersecting paths and the determinant of the distance matrix of a tree}

\date{\today}
\author[E. Briand]{Emmanuel Briand}
\address[E. Briand]{Departamento Matemática Aplicada 1, Universidad de  Sevilla}
\email{ebriand@us.es}
\author[L. Esquivias-Quintero]{Luis Esquivias-Quintero}
\address[L. Esquivias-Quintero]{Departamento de Álgebra, Universidad de Sevilla, and Departamento de {\'A}lgebra, Geometr{\'{\i}}a y Topolog{\'{\i}}a, Universidad Complutense de Madrid.}
\email{lesquivias@us.es}
\author[Á. Gutiérrez]{Álvaro Gutiérrez}
\address[Á. Gutiérrez]{Departamento de Álgebra, Universidad de Sevilla, and School of Mathematics, University of Bristol}
\email{a.gutierrezcaceres@bristol.ac.uk}
\author[A. Lillo]{Adrián Lillo}
\author[M. Rosas]{Mercedes Rosas}
\address[A. Lillo, M. Rosas]{Departamento de Álgebra, Universidad de Sevilla}
\email[A. Lillo]{alillo@us.es}
\email[M. Rosas]{mrosas@us.es}

\begin{document}

\begin{abstract}
We present the first combinatorial proof of the  Graham--Pollak formula  for the determinant of the distance matrix of a tree, via sign-reversing involutions and the Lindström--Gessel--Viennot Lemma. 
Our approach provides a cohesive and unified framework for the understanding of the existing generalizations and
$q$-analogues of the Graham--Pollak formula, and facilitates the derivation of natural simultaneous generalizations for them.
 \end{abstract}

\maketitle

\section{Introduction}

Consider a tree $T$ with vertices labeled from $1$ to $n$ and edge set $E$.  The {\em distance} $d(i,j)$ between two vertices $i$ and $j$ is defined as the number of edges in the unique path connecting them in $T$. The {\em distance matrix} of $T$ is then defined as $M(T) = \left(d(i,j)\right)_{1 \le i, j \le n}$.
In their influential 1971 paper \cite{GrahamPollak}, Graham and Pollak  established that the determinant of the distance matrix of $T$ obeys what is now known as the {\em Graham--Pollak formula:}
\[
    \det M(T) = (-1)^{n-1}(n-1)2^{n-2}.
\]
Observe that the formula is solely dependent on the number of vertices of the tree, and not on the tree structure.

Multiple techniques from linear algebra, ranging from Gaussian elimination to Charles Dodgson's condensation formula, have been used to prove the Graham--Pollak formula \cite{GrahamPollak, YanYeh:2007, ZhouDing, DuYeh}. The expression $(n-1)2^{n-2}$ strongly suggests that $\det M(T)$ enumerates something; however, none of these proofs are combinatorial.
The formula has also been subject to various deformations and generalizations which we review below. 

The three main results of this present work are (i) the first combinatorial proof of the Graham--Pollak formula, (ii) a new highly general deformation of the formula, and (iii) a new highly general deformation of the formula that still depends on nothing but the number of vertices of the tree.\medskip

We interpret the distance between vertices of $T$ as the number of marked paths (paths  with a distinguished step) between them, turning the evaluation of the determinant into a signed enumeration problem for families of marked paths of $T$.  
We aim to utilize the powerful and sign-reversing Lindström--Gessel--Viennot (LGV) involution, which is only applicable when the problem is framed as a signed enumeration of families of ordinary paths. Therefore, we search for a way to associate to  $T$ a  network  such that: (i) the marked paths of $T$ lift to ordinary paths in this network; (ii) the network is suitable for applying the Lindström--Gessel--Viennot Lemma;  (iii) the involution is strong enough so that its fixed points have all the same sign. 

In order to fulfill (i), we build the network in two levels, {\em South} and {\em North}, so that the tail of any marked path is lifted to a path in  South and its head to a path in North.
For (ii) and (iii),  we blow up each vertex of the tree in each level into an ear, from which its adjacent edges hang like earrings. See Figure \ref{fig: route map full} for an illustration.

In order to fulfill (iii), it turns out one has to first split the families of marked paths into classes, and define a different network for each class.
  
After setting aside classes with simpler zero-sum involutions, the LGV involution leaves exactly one surviving element per class, and all the surviving elements share the same sign.

This paper is structured as follows: after some basic definitions in Section \ref{sec: basis},
our journey begins in Section \ref{sec: catalysts}, where we introduce {\em catalysts} for a tree and interpret $\det M(T)$ as the signed enumeration of all catalysts for $T$. We define the {\em arrowflow} induced by a catalyst in Section \ref{sec: arrowflows}. We partition the set of catalysts according to their induced arrowflows, which are classified to be either 
{\em zero-sum} or {\em unital}.
Our first contribution is Theorem \ref{thm: sum on class}.
\begin{mainthm}\label{thm: sum on class}
    Let $A$ be an arrowflow, let $C(A)$ be the set of catalysts for $T$ inducing $A$. Then
    \[
    \sum_{\kappa \in C(A)} \sign(\kappa) = \begin{cases}
        (-1)^{n-1} & \text{ if $A$ is unital,}\\
        0 & \text{ if $A$ is zero-sum.}
    \end{cases}
    \]
\end{mainthm}
We finish Section \ref{sec: arrowflows} by deriving the Graham--Pollak formula from Theorem \ref{thm: sum on class}, the proof of which constitutes the bulk of this work and unravels in the following sections.

The first part of Theorem \ref{thm: sum on class} is proved in Section \ref{sec: zero-sum} by means of a sign-reversing involution without fixed points on each zero-sum arrowflow class.

With the goal of proving the second part of Theorem \ref{thm: sum on class}, we reformulate our problem as a signed enumeration of path systems in a network. Given a unital arrowflow $A$, we construct an acyclic network $\route$ in Section \ref{sec: standard}, called the \emph{route network} of $A$.
In Sections \ref{subsec: lift and proj} and \ref{subsec: Lambda lifts npaths}, we establish a sign-preserving injection that sends catalysts with induced arrowflow $A$ into path systems on $\route$. 
Our injection satisfies the conditions required to  apply the Lindström--Gessel--Viennot Lemma \ref{lem: LGV}. Thus, the signed sum of all catalysts in the arrowflow class of $A$ is equal to the signed count of non-intersecting path systems on $\route$. 
In Sections \ref{NonintersectingHemisphere} and \ref{NonintersectingRoute}, we discuss existence and uniqueness of a non-intersecting path system on $\route$.
The unique non-intersecting path system for a given unital arrowflow is retrieved by a {\em Depth-First-Search walk} on $T$. From this, it follows that its underlying permutation is an $n$-cycle. 
This argument concludes our proof of Theorem \ref{thm: sum on class}. The details of some technical proofs are postponed to Section \ref{technical proofs}.

Section \ref{sec: generalizations} is devoted to deformations and generalizations of the distance matrix and of  the Graham--Pollak Formula. 
Replacing the distances $d(i,j)$ by their $q$-analogues $1+q+q^2+\cdots +q^{d(i,j)-1}$, one gets a matrix whose determinant is given by
$
(-1)^{n-1}(n-1)(1+q)^{n-2}
 $~
\cite{BapatLalPati, YanYeh:2007}. 
The determinant is, again and remarkably, independent on the structure of the tree.
Other deformations with the same property are obtained by putting weights on the edges \cite{BapatKirklandNeumann} or on the arcs \cite{BapatLalPati2009, ZhouDing}. 
Combining the two deformations, one obtains weighted $q$-analogues with 
weights on the edges \cite{YanYeh:2007, BapatLalPati} or on the arcs \cite{LSZ}.
More recently, Choudhury and Khare found a very general formula \cite[Thm.~A and Rem.~1.10]{CK19} that specializes to all of the above.

Our approach to the Graham--Pollak formula is particularly well-suited for studying deformations, which turn into weighted enumeration of combinatorial objects.
We interpret the distance $d(i,j)$ as the number of \emph{marked paths} (paths with a distinguished step) from $i$ to $j$.
The {\em weight of the marked path} 
\begin{equation}\label{eq: marked path}
i_0~i_1\ldots i_{k-1}~\underline{i_k~i_{k+1}}~i_{k+2} \ldots i_{d-1}~i_d
\end{equation}
is defined as the monomial
$x_{i_0i_1}\cdots x_{i_{k-1}i_{k}}y_{i_ki_{k+1}}z_{i_{k+1}i_{k+2}}\cdots z_{i_{d-1}i_{d}}$  
in three families of commuting variables attached to the arcs of the tree.  
The marked distance $d'(i,j)$ between $i$ and $j$ is defined as the sum of the weights of,
all marked paths from $i$ to $j$. 
Our second contribution is Theorem \ref{main: formula Emmanuel}, a compact formula for the determinant of the matrix $M'(T)$ of the marked distances.
\begin{mainthm}\label{main: formula Emmanuel}
Under the hypothesis 
$x_{ji} = 1/x_{ij}$ for all edges $\{i, j\}$ of $T$, 
the determinant of marked distance  matrix of $T$ is
\[
(-1)^{n-1} \sum_{e = \{a, b\}\in E} y_{ab}y_{ba} \prod_{(i,j) \in U(e)} (y_{ij} x_{ji} + y_{ji} z_{ij}),
\]
where $U(e)$ is the set of arcs supported on $T$ and ``pointing to $e$''.
\end{mainthm}

Notably, $\det M'(T)$ is no longer independent of the tree structure of $T$. 
By imposing certain simple relationships among the variables, we restore the independence of the determinant $\det M'(T)$ from the tree structure of $T$. This is our third contribution, Theorem \ref{main: formula indep}, which we will shortly state.

Towards our goal, the weight of the marked path \eqref{eq: marked path} is defined as
\begin{equation}\label{eq: simple metric}
x_{i_0i_1}\cdots x_{i_{k-1}i_{k}}\alpha_{\{i_k,i_{k+1}\}}(z_{i_k i_{k+1}} - x_{i_k i_{k+1}})
z_{i_{k+1}i_{k+2}}\cdots z_{i_{d-1}i_{d}}.
\end{equation}
where we attach to each edge $e$ one variable $\alpha_e$, and to each arc $\gamma$ two variables $x_\gamma$ and $z_{\gamma}$.
For each edge $e$ let $e^+$ and $e^-$ denote the two arcs (orientations) supported on $e$, and impose $x_{e^-}=1/x_{e^+}$.
We define the {\em generalized distance matrix} $M_G(T)$ as the matrix whose  $(i,j)$ entry is the sum of all the weights \eqref{eq: simple metric} of the marked paths from $i$ to $j$.

\begin{mainthm}\label{main: formula indep}
The determinant of the generalized distance matrix $M_G(T)$, specialized at $x_{e^+}=x_e$ and $x_{e^-}=1/x_e$, is equal to
\begin{align*}
(-1)^{n-1} 
\sum_{e\in E} \alpha_{e}^2  \big(z_{e^+}-x_e\big) \big(z_{e^-} - x_e^{-1}\big)
\prod_{\substack{f\in E\\f \neq e}}
\alpha_{f} (z_{f^+} z_{f^-} - 1).
\end{align*}
where  $e^+$ and $e^-$ denote the two arcs supported on  $e$. 

\end{mainthm}
As a result, note that even if $M_G(T)$ depends on the structure of the tree, on the choice of the arcs $e^+$ and $e^-$ for each edge $e$, and on the assignments of the weights to the edges, its determinant does not depend on any of these. In particular, it is invariant under all permutations of the families of weights $(\alpha_e, x_e, z_{e^+}, z_{e^-})$ assigned to each edge.

Finally, observe that the Choudhury--Khare deformation \cite[Thm.~A, case $x=0$]{CK19} can be recovered by specializing \( x_e \) to \( 1 \) in \( M_G(T) \). Consequently, all specializations of the Choudhury--Khare deformation can also be obtained from Theorem \ref{main: formula indep}.

\section{Preliminaries}
 \label{sec: basis}

\subsection{Basic definitions}

We assume familiarity with the standard concepts of graph theory. Nonetheless, for the sake of clarity, we provide definitions for some of the terms utilized in this work. 

We say \emph{edge} for non-oriented edge, \emph{arc} for oriented edge, \emph{graph} for simple graph and \emph{digraph} for directed graph.  
In a digraph, a vertex $i$ is a \emph{predecessor} (resp.~\emph{successor}) of vertex $j$ if $\arc{i}{j}$ (resp.~$\arc{j}{i}$) is an arc of the digraph. 
The \emph{arcs supported on an edge} $\{i,j\}$ of a simple graph $G$ are the ordered pairs $\arc{i}{j}$ and $\arc{j}{i}$. We denote with $\arcs{G}$ the set of all \emph{arcs supported on} $G$, i.e.~the arcs supported on the edges of $G$.
We consider \emph{vertex} and \emph{node} as synonyms, but we try to consistently use \emph{vertex} for the trees considered in this work ($T$, $T_A$, $Y$, $\earring{Y}$), and \emph{node} for the digraphs ($\hemisphere$, $\north$, $\south$, $\route$) we build from them, since some vertices (nodes) of the latter represent edges of the former.

A  {\em walk} on a graph (resp.~on a digraph) is a sequence of vertices where each pair of consecutive nodes is an edge (resp.~an arc).
Given a walk $x=x_0 x_1 \cdots x_m$, we say that node $x_0$ is its  {\em origin}, node $x_m$ its  {\em terminus}, and that the ordered pairs $\arc{x_i}{ x_{i+1}}$ are the  {\em steps} of the walk.
Finally, a  {\em path} is a walk with no repeated vertices.  
Given two vertices $i$ and $j$ of a tree $U$, we write  $\path{i}{j}{U}$ to denote the unique path of $U$ from $i$ to $j$.

We define a marked path in a  graph as a pair $(\pi; \gamma)$ where $\pi$ is a path and $\gamma$ is a step of $\pi$, that we call the \emph{marked step}. Writing $\pi=j_0j_1\cdots j_m$, there exists $k <m$ such that $\gamma=\arc{j_k}{j_{k+1}}$. 
We refer to the subwalks
$j_0 j_1 \cdots j_k$ and $j_{k+1} j_{k+2}\cdots j_m$
that  decompose $\pi$ as the  {\em tail} and  {\em head} of $(\pi;\gamma)$.  The marked step $\arc{j_k}{j_{k+1}}$ does not belong neither to the tail nor to the head.

A {\em network} $N$ is a digraph with two distinguished sequences $\Delta$ and $\nabla$ with an equal  number of (distinct) nodes: the  {\em sources} $(\Delta_1, \ldots, \Delta_n)$ and the  {\em sinks}
$(\nabla_1,\ldots,\nabla_n)$.
We require for no node to be simultaneously a source and a sink. We do not require that sources (resp. sinks) have only outgoing arcs (resp. incoming arcs) as incident arcs.

A {\em path system} of $N$ is a set of $n$ paths, where each path starts from a different source, and ends at a different sink.  
A path system induces a permutation $\sigma$, determined by the fact that for all $i$, the  path that starts at source $\Delta_i$ ends at sink $\nabla_{\sigma(i)}$. The {\em sign} of a path system is the sign of the corresponding permutation. A path system such that no node of the network belongs to more than one of its paths is called \emph{non-intersecting}.
We denote with $\PS(N)$ the set of the path systems of $N$, and with $\NIP(N)$ the subset of its non-intersecting path systems.

We denote by $\steps{\lambda}$ the set of all steps of a path $\lambda$. 
Given a path system $\Lambda = \{\Lambda_1, \ldots, \Lambda_n\}$ in a network,
we define $\steps{\Lambda}$ as the multiset of the steps of the $\Lambda_i$, where the multiplicity of a step is the number of paths to which it belongs.

\subsection{Sign-reversing involutions and the Lindström--Gessel--Viennot Lemma}

Given a finite set $X$ with a sign function $s:X \rightarrow \{-1, +1\}$, the \emph{sign enumeration problem} for $X$ consists in computing the sum  $\sum_{x \in X} s(x)$. A standard technique for solving it consists in exhibiting a \emph{sign reversing involution}. This is an involution $\phi$ on $X$ that sends any non-fixed point to an element of the opposite sign. Terms of the sum paired  by $\phi$   cancel out  and the sign enumeration problem simplifies as follows:
\[
\sum_{x \in X} s(x) = \sum_{x \in F} s(x)
\]
with $F$ the set of fixed points of $\phi$. If the sign reversing involution is strong enough, all $x$ in $F$ have the same sign and the sum further reduces to $\pm \# F$, that is, an unsigned enumeration problem.

Given a network $N$,  there is a sign-reversing involution on $\PS(N)$ that consists (skipping details) in swapping the tails (or the heads) of two of the paths that intersect, if any. Its fixed points are the non-intersecting path systems of $N$. Lindström and Gessel and Viennot \cite{Lindstrom, Gessel-Viennot} have applied widely this involution, for the enumeration of  combinatorial determinants.

We will make use of the following properties of this involution. 

\begin{lem}[Lindström--Gessel--Viennot Involution] 
\label{lem: LGV}
For any  acyclic network $N$,
there exists a sign-reversing involution on the set of its path systems  with the following properties: 
\begin{myenum}
\item Its fixed points are the non-intersecting path systems.
\item It preserves the multiset of steps of the path systems.
\end{myenum}
\end{lem}
We refer informally to this involution as the ``LGV involution'' of the network $N$.

As a consequence of the existence of the LGV involution, 
\begin{equation*}\label{eq: LGV simplification}
\sum_{\Lambda \in \PS(N)} \sign(\Lambda)
=
\sum_{\Lambda \in \NIP(N)} \sign(\Lambda).
\end{equation*}

The LGV involution can be applied as well to  networks whose arcs carry weights. 
Then the weight of a path  is defined as the product of the weights of its steps, and the weight $w(\Lambda)$ of a path system $\Lambda$ as the product of the weights of its paths.
The lemma now implies
\begin{equation*}\label{eq: LGV simplification with weights}
\sum_{\Lambda \in \PS(N)} \sign(\Lambda) w(\Lambda)
=
\sum_{\Lambda \in \NIP(N)} \sign(\Lambda)
w(\Lambda).
\end{equation*}


\section{Catalysts}\label{sec: catalysts}

Throughout the paper, we fix a tree $T$ with vertex set $V(T)=\{1, 2, \ldots, n\}$ 
and edge set $E$. 

Given a permutation $\sigma$ of $V(T)$ and a map $f : V(T) \rightarrow \arcs{T}$,  we say that the ordered pair  $(\sigma, f)$ is a   {\em catalyst for $T$} if for each vertex $i$, its image $f(i)$ is a step of the path $\pathT{i}{\sigma(i)}$.
The {\em sign} of a catalyst is the sign of its underlying permutation $\sigma$. 
Catalysts  encode the sets  of $n$ marked paths $\left(\pathT{i}{\sigma(i)}; f(i)\right)$ whose origins, as well as terminus, are all $n$ vertices of $T$. 

\begin{ex}\label{example_signed_path}
Consider the tree $T$ of Figure \ref{fig: catalyst T}. Figures \ref{fig: catalyst A}, \ref{fig: catalyst B}, \ref{fig: catalyst C} show the sets of marked paths corresponding to three catalysts for $T$. 

\begin{figure}[h]
        \begin{subfigure}[b]{0.25\textwidth}
        \begin{center}
        \includegraphics[width = \textwidth, page=1]{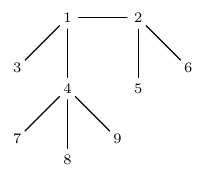}
        \end{center}
        \caption{}
        \label{fig: catalyst T}
        \end{subfigure}
        \begin{subfigure}[b]{0.24\textwidth}
        \includegraphics[width = \textwidth, page=2]{EarringGraph_all_pictures.pdf}
        \caption{}
        \label{fig: catalyst A}
        \end{subfigure}
        \begin{subfigure}[b]{0.24\textwidth}
        \includegraphics[width = \textwidth, page=3]{EarringGraph_all_pictures.pdf}
        \caption{}   
        \label{fig: catalyst B}
        \end{subfigure}
        \begin{subfigure}[b]{0.24\textwidth}
        \includegraphics[width = \textwidth, page=4]{EarringGraph_all_pictures.pdf}
        \caption{}
        \label{fig: catalyst C}
        \end{subfigure}
    \caption{A tree $T$ (\ref{fig: catalyst T})  and the sets of marked paths corresponding to three  of its catalysts (\ref{fig: catalyst A}, \ref{fig: catalyst B}, \ref{fig: catalyst C}).}
    \label{fig: catalysts}
    
\end{figure}
For instance, diagram \ref{fig: catalyst A} depicts the catalyst $(\sigma, f)$ where $\sigma$ and $f$ are given by the following table.
\[
\begin{array}{c|ccccccccc}
i & 1 & 2 & 3 & 4 & 5 & 6 & 7 & 8 & 9 \\\hline
\sigma(i) & 6 & 5 & 8 & 7 & 3 & 2 & 9 & 4 & 1 \\
f(i) & 12 & 25 & 31 & 47 & 21 & 62 & 49 & 84 & 41 
\end{array}
\]
\end{ex}

Now, since  $d(i,\sigma(i))$ counts the steps of the unique path of $T$ from $i$ to $\sigma(i)$, 
the determinant $\det M(T)$ does a signed enumeration of  all catalysts for $T$. 
Indeed, by definition
\begin{equation}
\det M(T) = \sum_{\sigma} \sign(\sigma) \, d(1,\sigma(1))d(2,\sigma(2)) \cdots d(n,\sigma(n)).\label{eq:P}   
\end{equation}
where the sum is over all permutations $\sigma$ of $V(T)$. 
Therefore,
$
\det M(T)
= \sum_\kappa\sign\kappa
\label{eq:prod}
$, 
where  the sum ranges over all catalysts for $T$.  

It is worth noting that the definition of catalyst implies that its underlying permutation  is a derangement (a permutation without  fixed points).

\section{Arrowflows and the Graham--Pollak formula}\label{sec: arrowflows}
  
\begin{figure}[h!]
    \begin{subfigure}{0.3\textwidth}
        \includegraphics[width = \textwidth, page=5]{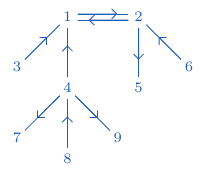}
        \caption{Unital\\\strut}
        \label{fig:arrowflow A}
\end{subfigure}
\hfill
\begin{subfigure}{0.3\textwidth}
        \includegraphics[width = \textwidth, page=6]{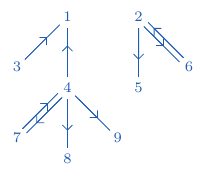}
      \caption{Disconnected\\ zero-sum}
      \label{fig:arrowflow B}
\end{subfigure}
\hfill
\begin{subfigure}{0.3\textwidth}
        \includegraphics[width = \textwidth, page=7]{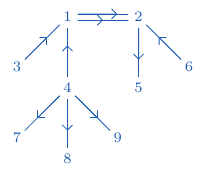}
      \caption{Connected\\ zero-sum}
      \label{fig:arrowflow C}
\end{subfigure}
\caption{Unital and zero-sum arrowflows.}\label{fig: arrowflows}
\end{figure}

  Let us define an {\em arrowflow on $T$} as a multiset of $n$ arcs supported on $T$. 
Each catalyst $\kappa = (\sigma, f)$  induces an arrowflow on $T$:
the multiset defined as the image of $f$, where the multiplicity of an arc $\gamma$ is the number of vertices $i$ such that $f(i) = \gamma$.
Remark that we draw arrowflows using $T$ as blueprint, as illustrated in Figure 
\ref{fig: arrowflows}.

Different catalysts on $T$ can result on the same arrowflow. On the other hand, there exist arrowflows on $T$ that are not induced by any catalyst for $T$. 
The set of catalysts inducing an arrowflow $A$ on $T$ defines  the {\em arrowflow class of  $A$}, denoted by $C(A)$.

 The {\em arrowflow partition} is the partition of the set of all catalysts for $T$ formed by the nonempty arrowflow classes. Rewriting the signed enumeration of catalysts done by $\det M(T)$ in \eqref{eq:P} according to the arrowflow partition, we obtain
\begin{align} \label{main_eq}
\det M(T) = \sum_{\substack{A\\ \text{arrowflow}} }\sum_{\kappa \in C(A)} \sign(\kappa),
\end{align}
where the first sum is taken over all arrowflows on $T$, and the second one over all catalysts $\kappa$ in the arrowflow class $C(A)$. Empty arrowflow classes have no effect on this summation.

It turns out to be convenient to introduce a  classification of arrowflows.
We say that an arrowflow $A$ is {\em connected} when each edge $\{i,j\}$ of $T$ supports an arc $\arc{i}{j}$ or $\arc{j}{i}$ in $A$. If $A$ is a connected arrowflow, then there exists precisely one edge of $T$ supporting two arcs of $A$. We call them  the {\em marked arrows} and the underlying edge the {\em marked edge}.
An arrowflow is said to be {\em unital} when it is connected and has no repeated arcs, as illustrated in Figure \ref{fig:arrowflow A}. Otherwise, it is said to be {\em zero-sum}. There are two possible causes for an arrowflow to be zero-sum: either the arrowflow is disconnected, as illustrated in 
 Figure \ref{fig:arrowflow B}, or
the arrowflow is connected, and there is a repeated arrow as illustrated in
Figure \ref{fig:arrowflow C}.

The arrowflow partition is the optimal way of partitioning the set of catalysts for $T$, as made precise in Theorem \ref{thm: sum on class}, which we restate below.

\setcounter{mainthm}{0}
\begin{mainthm} 
Let $A$ be an arrowflow. 
Then
\begin{align}\label{crucial_result}
\sum_{\kappa \in C(A)} \sign(\kappa) = 
\begin{cases}
(-1)^{n-1} & \text{ if } A \text{ is unital,}\\
0 & \text{ if } A  \text{ is zero-sum.}
\end{cases}
\end{align}
\end{mainthm}

The following sections will present a combinatorial proof of this result. We close this section by deriving the Graham--Pollak formula from it.

\begin{proof}[Derivation of the Graham--Pollak Formula from Theorem \ref{thm: sum on class}]

After \eqref{main_eq} and \eqref{crucial_result}, 
\[
\det M(T) = (-1)^{n-1} \times \text{number of unital arrowflows on $T$.} 
\]
The number of unital arrowflows on $T$ is $(n-1) \, 2^{n-2}$ :  the factor $(n-1)$ counts the  ways of selecting the marked edge, whereas the factor $2^{n-2}$ counts the number of ways in which  the remaining $n-2$ edges can be oriented. Whence,
\[
\det M(T) = (-1)^{n-1} (n-1) 2^{n-2}. \qedhere
\]

\end{proof}
\section{Zero-sum arrowflows}
\label{sec: zero-sum}

Let $A$ be a zero-sum arrowflow. 
We show that the signed sum of catalysts that induce $A$ is zero by constructing a sign-reversing involution $\varphi_A$ of $C(A)$ without fixed points. 
This amounts to one half of Theorem \ref{thm: sum on class}. The argument is split into two cases depending on whether $A$ is connected or not.

\begin{lem}\label{prop zero-sum disconnected}
Let $A$ be a disconnected zero-sum arrowflow on $T$.
Fix  an edge  $\{i,j\}$ of $T$ among those carrying no arc of $A$.  
For any $(\sigma, f) \in C(A)$, set 
$\varphi_A(\sigma,f)  =(\sigma \circ (i ~ j),f \circ (i ~j))$, where $(i ~j)$ denotes the transposition swapping $i$ and $j$.

Then $\varphi_A$ is a sign-reversing involution of $C(A)$ without fixed points.
\end{lem}

\begin{figure}[h]
    \centering
        \includegraphics[height = 3cm, page=8]{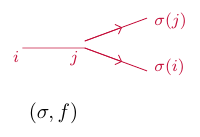}
        \hspace{2em}
        \includegraphics[height = 3cm, page=9]{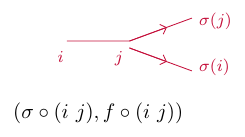}
    \caption{Involution $\varphi_A$  on a zero-sum disconnected arrowflow  in which $\{i,j\}$ does not separate $\sigma(i)$ and $\sigma(j)$.}
    \label{fig: disconnected 1}
\end{figure}
\begin{figure}[h]
    \centering
        \includegraphics[height = 3cm, page=10]{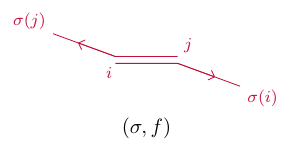}
        \hspace{2em}
        \includegraphics[height = 3cm, page=11]{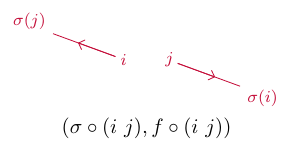}
    \caption{Involution $\varphi_A$ on a zero-sum disconnected arrowflow in which $\{i,j\}$ separates $\sigma(i)$ and $\sigma(j)$.
    }
    \label{fig: disconnected 2}
\end{figure}
\begin{proof}
It is enough to show that, for every $\kappa =(\sigma,f) \in C(A)$,  $\varphi_A(\kappa)$ is a catalyst. 
Indeed, since the multiset of the arcs $f(i)$, for $i \in V(T)$, 
is invariant under $\varphi_A$, we will conclude that $\varphi_A(\kappa) \in C(A)$. 

Let $\kappa=(\sigma,f) \in C(A)$ and $\varphi_A(\kappa)=(\tau,g)$.
By construction, for any vertex $k$ of $T$ distinct from $i$ and $j$, the arc $g(k)$ is a step in $\pathT{k}{\tau(k)}$.
It remains to show that it is also the case when $k$ is $i$ or $j$. We will prove it only for $k=i$, since the case $k=j$ is similar.

Let $\arc{a}{b}=f(j)$. Let $T'$ the tree obtained from $T$ by deleting the edge $\{a, b\}$.
For any vertex $k$, observe that $\arc{a}{b}$ is a step in $\pathT{k}{\sigma(k)}$ if and only if  $k$ and $a$ are in one connected component of $T'$ while $\sigma(k)$ and $b$ in the other. 
Since $\arc{a}{b}=f(j)$, it  is a  step of $\pathT{j}{\sigma(j)}$. 
Therefore 
$j$ and $a$ are in one component of $T'$, and $\sigma(j)$ and $b$ are in the other. 
Since $\{i, j\} \neq \{a, b\}$ (because $\{i, j\}$ supports no arc of $A$ but $\{a, b\}$ supports $\arc{a}{ b}=f(j)$),
the vertex $i$ is still a neighbor of $j$ in $T'$.
Therefore it is in the same connected component as $j$ and $a$. 
Applying the property again, we obtain that $\arc{a}{b}$, which is $g(i)$, is a step in $\pathT{i}{\sigma(j)}$, which is $\pathT{i}{\tau(i)}$.
\end{proof}

We remark that the involution described above does not change the multiset of steps of the catalyst whenever $\{i,j\}$ is \emph{not} in the path $\pathT{\sigma(i)}{\sigma(j)}$, see Figure \ref{fig: disconnected 1}. However, this is no longer true when we drop this assumption, see Figure \ref{fig: disconnected 2}. 

For connected zero-sum arrowflows, we introduce a similar involution. 
\begin{lem}\label{prop zero-sum connected}
Let $A$ be a connected zero-sum arrowflow on $T$.
Let $\arc{a}{ b}$ be the unique element of $A$ with multiplicity $2$.
For any $(\sigma, f) \in C(A)$, set 
$
\varphi_A(\sigma, f)
=\left(\sigma \circ (i ~ j),f\right),
$
where   $i$ and $j$ are the two vertices mapped by $f$ to  $\arc{a}{b}$.

Then $\varphi_A$ is a sign-reversing involution of $C(A)$ without fixed points.
\end{lem}
\begin{figure}[h]
    \includegraphics[height = 2.9cm, page=12]{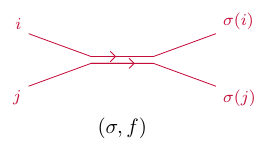}
    \hspace{2cm}
    \includegraphics[height = 2.9cm, page=13]{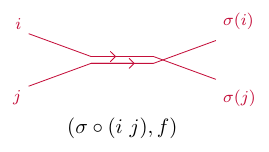}
    \caption{Involution $\varphi_A$ on a zero-sum connected arrowflow.}\label{fig: double edge}
\end{figure}
\begin{proof}
As in the disconnected case, it is sufficient to show that $f(i)$ is a step in $\pathT{i}{\sigma(j)}$. By symmetry, we also have that  $f(j)$ is a step in $\pathT{j}{\sigma(i)}$.
We proceed in the same way as the previous lemma: since $\arc{a}{b}=f(i)=f(j)$  is a step in both $\pathT{i}{\sigma(i)}$ and $\pathT{j}{\sigma(j)}$,  vertices $i,j$ and $a$ belong to one connected component of the tree obtained from $T$ by deleting $\{a,b\}$, whereas $\sigma(i)$, $\sigma(j)$ and $b$ are in the other one. See Figure \ref{fig: double edge}. 
\end{proof}
\section{The route network of an unital arrowflow}\label{sec: standard}

The proof of the second part of Theorem \ref{thm: sum on class}, which addresses unital arrowflows, is more nuanced, and  will be presented in the next few sections. 
In this section, we introduce an acyclic network $\route$ that we call the {\em route network} of the unital arrowflow $A$ on $T$.

To construct the route network, we first turn $T$ into a \emph{left-right tree} $T_A$ (a generalized binary tree where each vertex may have several left children and several right children).
From $T_A$ and its mirror image $T'_A$, we apply successively two basic constructions (\emph{earring tree} and \emph{restricted line digraph}) to get the \emph{hemisphere digraphs}: the \emph{Southern hemisphere} and the \emph{Northern hemisphere}. 
The route network is constructed on the union of these hemispheres, connected by arcs we call \emph{bridges}. Bridges correspond to the arcs in the arrowflow $A$, and always go from the Southern to the Northern hemisphere.

 In this section, we fix a unital arrowflow $A$ of $T$ with marked edge $\{a, b\}$.

\subsection{The rooted tree $T_A$}\label{sec: root insertion} 
Using our unital arrowflow $A$ we transform $T$ into a rooted tree $T_A$.
This is done by adding a new vertex $r$ and replacing the marked edge $\{a, b\}$ of $T$ with the pair of edges $\{r, a\}$ and $\{r, b\}$.   $T_A$ is rooted by declaring $r$ to be its root.

\begin{figure}[h!]
\begin{subfigure}[t]{0.33\textwidth}
\begin{center}
        \includegraphics[width = .8\textwidth, page=14]{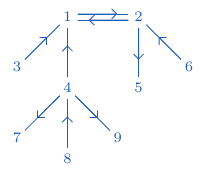}
\end{center}
\caption{\smaller Unital~arrowflow~$A$~with marked~edge~$\{1,2\}$~and marked arrows $12$ and~$21$.}
\label{fig: plane rooted tree A}
\end{subfigure}
\hfill
\begin{subfigure}[t]{0.3\textwidth}  
        \includegraphics[width = \textwidth, page=15]{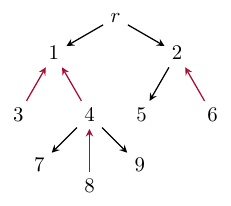}
\caption{\smaller  The rooted tree $T_A$ with its orientation $A_0$.}
\label{fig: plane rooted tree B}
\end{subfigure}
\hfill
\begin{subfigure}[t]{0.3\textwidth}   
        \includegraphics[width = \textwidth, page=16]{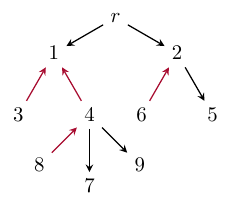}
\caption{\smaller  The left-right tree $T_A$.
}
\label{fig: plane rooted tree C}
\end{subfigure}
\caption{}
    \label{fig: plane rooted tree}\label{Running_example1}\label{T0}
\end{figure}

After setting $\rootins(\arcInput{a}{b})=\arc{r}{a}$, $\rootins(\arcInput{b}{a})=\arc{r}{ b}$ and $\rootins(\arcInput{i}{j})=\arc{i}{j}$ for any other  $\arc{i}{j}\in A$, we obtain from  arrowflow $A$ a set of arcs $A_0=\rootins(A)$. 
Since every edge $\{i, j\}$ in $T_A$ supports exactly one arc  of  $A_0$ (either $\arc{i}{j}$ or $\arc{j}{i}$), $A_0$ is an orientation of $T_A$.
See Figures \ref{fig: plane rooted tree A} and \ref{fig: plane rooted tree B}.

\begin{ex}\label{ex: running}
 Figure \ref{fig: plane rooted tree A}  shows  an arrowflow $A$ on tree $T$. Figure \ref{fig: plane rooted tree B} shows the corresponding rooted tree with its orientation $A_0$ induced by $A$. 
\end{ex}

Moreover, each marked path $(\pathT{i}{j}; \gamma)$ of $T$ is lifted into a marked path $(\pathTO{i}{j}; \rootins(\gamma))$ of $T_A$. 
The path $\pathTO{i}{j}$ is obtained from $\pathT{i}{j}$ by replacing any subword $ab$ with $arb$, and any  subword $ba$ with $bra$.

This lifting of marked paths induces a lifting of the catalysts in $C(A)$ (interpreted as sets of $n$ marked paths of $T$) into sets of $n$ marked paths of $T_A$.
Explicitly, the lifting to $T_A$ of $(\sigma;f)$ is the set 
\[
\left\{\left(\pathTO{i}{\sigma(i)}; \theta \circ f(i)\right) \text{ with $i \in V(T)$}\right\}.
\]

\subsection{The left-right tree of a unital arrowflow}\label{ConstructionT0}

A \emph{rooted plane tree} is a rooted tree in which the children of each vertex are totally ordered.
A {\em left-right tree} is a rooted plane tree where the set of children of each vertex is partitioned into two subsets: \emph{left children} and \emph{right children}, and where left children precede right children. 
This definition generalizes that  of a binary tree.  
The \emph{orientation} $\OriY$ of a left-right tree $Y$ is the set of arcs $\arc{i}{j}$ supported on $Y$ where either $i$ is a left child of $j$, or $j$ is a right child of $i$.

We endow $T_A$ with the structure of a left-right tree. For this, we first turn $T_A$ into a plane tree by choosing, for each vertex $i$, a total order $<_i$ on its children satisfying that $j <_i k$ whenever $\arc{j}{ i}$ and $\arc{i}{ k}$ are both in $A_0$. That is, up arrows should appear to the left of any down arrow.

 This makes of  $T_A$  a left-right tree with orientation $A_0$. A child $j$ of $i$ is a right child if and only if $\arc{i}{j}$ is in $A_0$. 

\begin{ex}\label{ex: running LR}
Figure \ref{fig: plane rooted tree C} shows the left-right tree for the rooted tree of Figure \ref{fig: plane rooted tree B}.
\end{ex}

\subsection{The earring tree of a left-right tree. }
\begin{figure}[h]
    \centering
    \begin{subfigure}[t]{.3\textwidth}
        \includegraphics[width = \textwidth, page = 17]{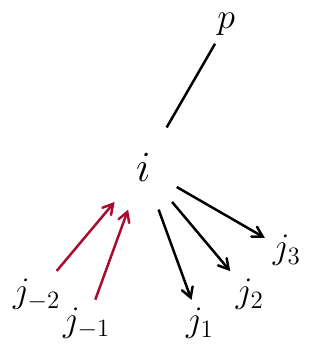}
        \caption{\;}\label{earring a}
    \end{subfigure}
    \begin{subfigure}[t]{.3\textwidth}
        \includegraphics[width = \textwidth, page = 18]{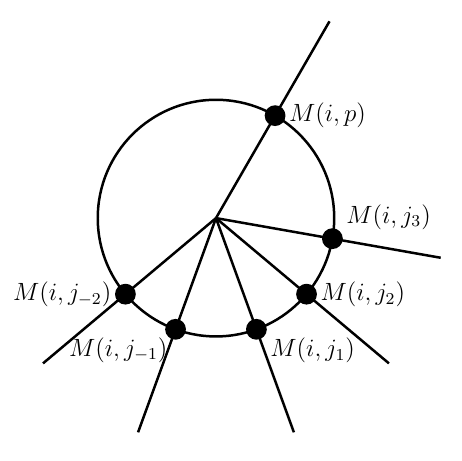}
        \caption{\;}\label{earring b}
    \end{subfigure}
    \begin{subfigure}[t]{.3\textwidth}
        \includegraphics[width = \textwidth, page = 19]{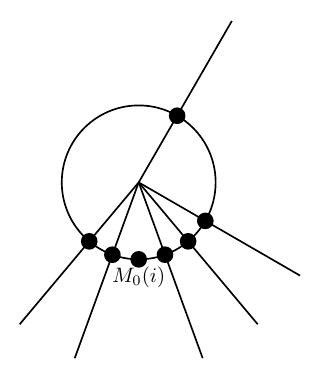}
        \caption{\;}\label{earring c}
    \end{subfigure}
    \begin{subfigure}[t]{.3\textwidth}
        \includegraphics[width = \textwidth, page = 20]{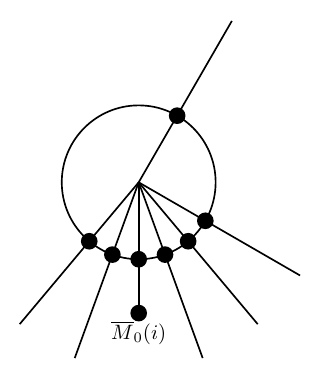}
        \caption{\;}\label{earring d}
    \end{subfigure}
    \begin{subfigure}[t]{.3\textwidth}
        \includegraphics[width = \textwidth, page = 21]{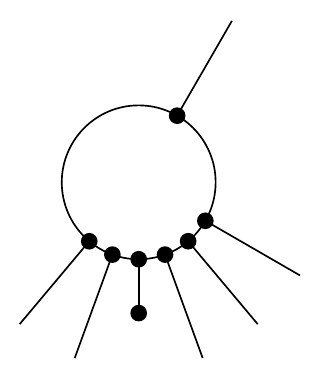}
        \caption{\;}\label{earring e}
    \end{subfigure}
    \begin{subfigure}[t]{.3\textwidth}
        \includegraphics[width = \textwidth, page = 22]{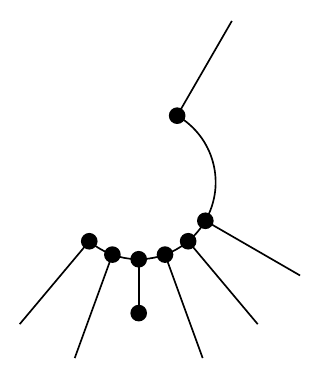}
        \caption{\;}\label{earring f}
    \end{subfigure}
    \caption{Steps in the construction of the earring tree.}
    \label{fig: local earring construction}
\end{figure}

The {\em earring tree} $\earring{Y}$ of a left-right tree $Y$ is the plane tree  obtained by ``blowing up'' the vertices of $Y$ different from the root as illustrated in Figure \ref{fig: local earring construction}.  
In the earring tree, each vertex $i$ of $Y$ is replaced with several nodes $M(i,j)$ (one for each neighbor $j$ of $i$)  plus two other vertices $M_0(i)$ and $\M_0(i)$.

Consider an embedding of  $Y$  in the plane in a way that preserves the plane tree structure. 
As one rotates counterclockwise around a node $i$, the left children appear first, followed by the right children, and finally the parent. See Figure \ref{fig: local earring construction}.

\eqref{earring a} Let $i$ be a vertex of $Y$ different from the root.
Let $p$ be its parent. 
Let
  \[
  j_{-\ell} <_i \cdots <_i j_{-2} <_i j_{-1} <_i j_1 <_i j_2 <_i \cdots <_i j_{m} 
  \]
  be the  children of $i$, labeled in such a way that   
   \( j_k \) is a left child of $i$ if \( k < 0 \), and a right child if \( k > 0 \).  

\eqref{earring b} Draw a circle of small radius $\varepsilon$ centered at $i$. Let $M(i,j_k)$ be the intersection of the circle with the edge $\{i, j_k\}$, and let $M(i,p)$ be the intersection of the circle with the edge $\{i, p\}$.

  \eqref{earring c} Insert a point $M_0(i)$ on the circle so that, when turning anticlockwise from $M(i,p)$, one finds first the $M(i,j)$ for $j$ the left children of $i$, then $M_0(i)$, and, finally, the $M(i,j)$ for $j$ the right children of $i$.

  \eqref{earring d} 
  Add a node $\M_0(i)$  symmetric of $i$ with respect to $M(i)$, and draw the segment joining $i$ to $\M_0(i)$.

\eqref{earring e}  Delete everything that is in the interior of the circle (the point $i$ and the segments joining it to points of the circle).

\eqref{earring f} The points $M(i,j_k)$ subdivide the circle into arcs.  Remove the first arc met by walking counterclockwise from $M(i,p)$.

\begin{figure}[h]
\hfill
\begin{subfigure}[t]{0.49\textwidth}
\begin{center}
        \includegraphics[width = .6\textwidth, page=23]{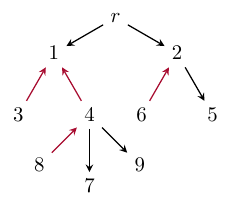}
\end{center}
\caption{A left-right tree $Y$}
\end{subfigure}
\begin{subfigure}[t]{0.50\textwidth}
\begin{center}
\includegraphics[height = 4cm, page=24]{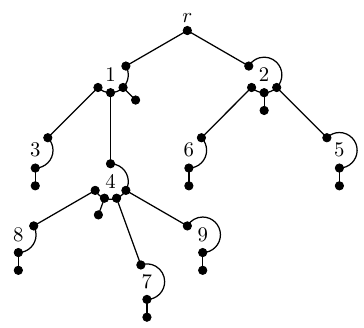}
\end{center}
\caption{The earring tree $\earring{Y}$}
\end{subfigure}
\caption{A left-right tree and its earring tree. }\label{fig: earring}
\label{Example_earring_tree}
\end{figure}

As a result of performing these blow-ups at all vertices of $Y$ different from the root, we obtained a plane tree embedded in the plane, whose vertices are $r$ and the $M_0(i)$, $\M_0(i)$ and $M(i,j)$. This is  the earring tree  $\earring{Y}$ of $Y$. 
See Figure \ref{fig: earring} for an example of the construction.

\subsection{The restricted line digraph}

We define the  \emph{restricted line digraph} $\linegraph{Z}$ of a graph $Z$ as the digraph whose nodes are the arcs $\arc{P}{ Q}$ supported on $Z$, and whose arcs are the pairs $(\arc{P}{ Q}, \arc{Q}{ R})$ where $P$, $Q$, $R$ are vertices of $Z$,  with $P$ and $R$ neighbors of $Q$ and $P \neq R$.  The condition $P \neq R$ does not appear in the standard definition of line digraph, hence the term \emph{restricted}. 

\begin{lem}[node-arc correspondence]\label{lem: proj from linegraph}
 If $Z$ is a tree, 
the map 
\[
(\arc{P_1}{ Q_1}\,,~ \arc{P_2}{ Q_2}\,, ~\ldots,~ \arc{P_d}{ Q_d}) \longmapsto (P_1, Q_1, Q_2, \ldots, Q_d)
\]
restricts to a bijection  from the set of all walks in $\linegraph{Z}$ to the set of all non-trivial paths in $Z$. 
\end{lem}

\begin{proof}
The walks in $\linegraph{Z}$ are  the sequences of arcs of $Z$ of the form
\[
\arc{P_1}{ P_2}\,,~ \arc{P_2}{ P_3}\,,~ \ldots, ~\arc{P_{m-1}}{ P_m}
\]
where, by the definition of the restricted line graph,   $P_{k+1}\neq P_{k-1}$ for all $k$.
The above map send them bijectively to the non-trivial walks with no recoil of $Z$. But since $Z$ is a tree, its walks with no recoil are its paths.
\end{proof}

\begin{cor}\label{cor unique walk}
  If $Z$ is a tree then $\linegraph{Z}$ has the  ``unique walk property'': 
between any two nodes of $\linegraph{Z}$, there is at most one walk. In particular, any walk of  $\linegraph{Z}$ is a path, and $\linegraph{Z}$ is acyclic.
\end{cor}

\subsection{The hemisphere network}

Given a left-right tree $Y$, the hemisphere network  $\hemi{Y}$ of $Y$ is defined as
\[
\hemi{Y}=\linegraph{\earring{Y}}.
\]
Let us denote with $V^*(Y)$ the set of all vertices of $Y$ different from $r$.
For each $i \in \pV(Y)$,  we set $v(i)=\M(i)M(i)$.
For two adjacent vertices $i$ and $j$ of $Y$, we will also write $e(\arcInput{i}{j})$ for $M(i,j)M(j,i)$.

We turn the hemisphere $\hemisphere$ into a a network by choosing the $v(i)$ as its sources, and the $e(\gamma)$ for $\gamma\in \OriY$ (the orientation set of $Y$) as its sinks. 
As a consequence of Corollary \ref{cor unique walk}, we have:
\begin{p}\label{south is acyclic} Let $\treeU$ be a left-right tree.
Given any pair of nodes $x$ and $y$ of $\hemi{\treeU}$, there is at most one walk in $\hemi{\treeU}$ from $x$ to $y$. 
Therefore, any walk in $\hemi{\treeU}$ is a path, and $\hemi{\treeU}$ is acyclic.
\end{p}

\subsection{The route network} \label{se:route}

We come back to the situation where we have a tree $T$, an arrowflow $A$, and the plane rooted tree $T_A$ with 
orientation $A_0$ obtained from them, as described in Section \ref{sec: root insertion}.
We now construct a digraph $\route$ out of two disjoint subgraphs $\south$ and $\north$ (Southern and Northern hemispheres), and a set of bridges connecting them, as we proceed to describe.

\begin{figure}[ht]
\begin{subfigure}[t]{0.49\textwidth}
\begin{center}
        \includegraphics[width = .6\textwidth, page=25]{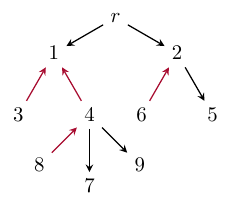}
\end{center}
\caption{\small The left-right $T_A$, oriented by $A_0$.}
\end{subfigure}
\hfill
\begin{subfigure}[t]{0.49\textwidth}
\begin{center}
        \includegraphics[width = .6\textwidth, page=26]{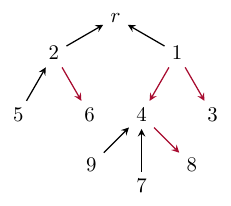}
\end{center}
\caption{\small The flipped image $T'_A$ of $T_A$, with orientation the arcs in $A_0$ reversed.}
\end{subfigure}
\caption{\small A left-right tree $T_A$ and its mirror image $T_A'$.}
\label{fig: mirror}
\end{figure}

The \emph{Southern hemisphere $\south$} is defined as 
\( 
\south = \hemi{T_A}.
\)
Similarly, the {\em Northern hemisphere $\north$} is defined as $\north=\hemi{T_A'}$, where $T_A'$ is a left-right tree $T'_A$, that can be thought out as the mirror image of $T_A$, see Figure \ref{fig: mirror}. 

 As a rooted tree, $T'_A$ is identical to $T_A$, but its is oriented in the opposite way. That is, $\Ori{(T_A')}$ is taken as the set of arcs in $A_0$ reversed.
Thus, the plane directed tree structure of $T_A$ is a flipped version of the one of $T_A$. Any right child of a vertex $ i $ becomes a left child in $ T'_A$, and vice-versa. 
The plane tree structure of $ T'_A$ is defined by saying that, for each vertex $ i $ in $ T'_A $, the order on the children of $ i $ is obtained by reversing the order $ <_i $ on the children of $ i $ in $ T_A $. It follows that $T'_A$ is also a left-right tree.

\begin{figure}[h]
    \centering
        \includegraphics[width = .8\textwidth, page=27]{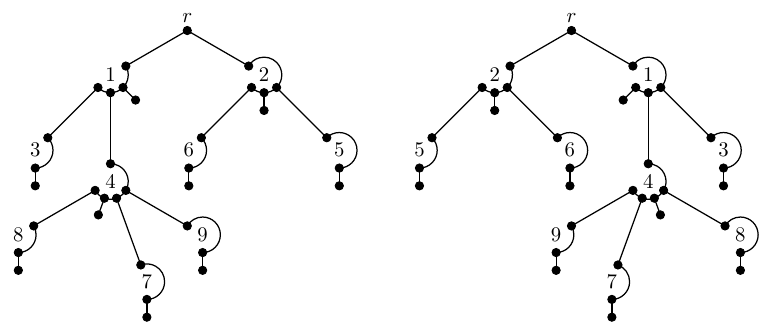}
    \caption{ The plane rooted trees $\earring{T_A}$ and $\earring{T_A'}$. The nodes of $\route$ are the arcs they support. }
    \label{fig: route map full}
\end{figure}

To build $\route$ we will consider the disjoint union $\south \sqcup \north$.  See Figure \ref{fig: route map full}. 
To distinguish between the homonymous nodes $v(i)$ and $e(\gamma)$ in $\south$ and $\north$,  we denote $v'(i)$ and $e'(\gamma)$ for $v(i)$  and $e(\gamma)$ when considered as nodes of $\north$, while we keep the notations $v(i)$ and $e(\gamma)$ for the nodes of $\south$.
Similarly, we prime the vertices of $\smash{\earring{T_A'}}$ to  distinguish them from those of $\earring{T_A}$. 
So $r'$, $M'_0(i)$, $\smash{\Mp_0}(i)$ and $M'(i,j)$ are the vertices of $\earring{T_A'}$ corresponding to the vertices $r$, $M_0(i)$, $\smash{\M}_0(i)$ and $M'(i,j)$ of $\earring{T_A}$.

The {\em route network $\route$} is obtained from $\south \sqcup \north$ by adding $n$ arcs $\arc{e(\gamma)}{ e'(\gamma})$, one for each $\gamma$ in $A_0$, connecting them. 
We call these arcs the {\em bridges} between hemispheres of  $\route$.

\begin{cor} \label{cor: paths in N and R} Any walk in either $\south$, $\north$ or $\route$ is  a path. 
In particular, $\south$, $\north$ and $\route$  are acyclic.
\end{cor}

\begin{proof}
This follows from 
Proposition \ref{south is acyclic} for $\south$ and $\north$.
Finally, since all bridges between hemispheres point from South to North,  any walk $w$ in $\route$ not contained in $\south$ or $\north$ factorizes uniquely as $w=w'w''$ with $w'$ a walk in $\south$ and $w''$ a walk in $\north$. Since $w'$ and $w''$ are paths, so is $w$.  
\end{proof}

\subsection{The network structures}

Recall that, for each $i\in V(T)$,  $v(i)=M_0(i)\M_0(i)$,  is a node of $\south$.  In analogy, we let $v''(i)=\Mp_0(i)M_0'(i)$,  it is a node of $\north$.
We make of $\route$ a network by taking sources $v(i)$ for $i \in V(T)$ and sinks $v''(i)$ for $i \in V(T)$.

Let $\Psi$ be the involution of $\north$ that sends every node $\arc{P}{Q}$  to $\arc{Q}{P}$ (remember that the nodes of $\north$ are arcs of $\earring{T'_A}$, so here $P$ and $Q$ are the vertices of $\earring{T'_A}$). 
For instance, $\Psi(v''(i))=\Psi(M_0'(i){\Mp_0}(i))={\Mp_0}(i) M_0'(i)=v'(i)$ and $\Psi(e'(\arcInput{i}{j}))=\Psi(M'(i,j)M'(j,i))=M'(j,i) M'(i,j)=e'(\arcInput{j}{i})$. 

The involution $\Psi$ is an anti-automorphism of $\north$. Given a walk $x_1 x_2\cdots x_n$ of $\north$, then $\Psi(x_n)\cdots \Psi(x_2)\Psi(x_1)$ is a walk of $\north$, which we denote it by $\Psi(x_1x_2\cdots x_n)$.

In Sections \ref{subsec: lift and proj} and \ref{subsec: Lambda lifts npaths}, we will lift  any marked path $(i\cdots j;pq)$ of $T$  with marked step $\arc{p}{q}$ belonging to $A$  to a path
\begin{equation}\label{eq: one lifting}
v(i) \cdots e(\arcInput{p}{q})e'(\arcInput{p}{q}) \cdots v''(j)
\end{equation}
of $\route$. This will turn the sign enumeration of catalysts in $C(A)$ into a sign enumeration of path systems in $\route$.

The lifting \eqref{eq: one lifting} factorizes as a concatenation $\mu \ \Psi(\mu')$, where $\mu$ is the path in the Southern hemisphere that starts in $v(i)$ and ends at $e(\arcInput{p}{q})$,
and $\mu'$ in the path of Northern hemisphere that starts at  $v''(j)$ and ends at $  e'(\arcInput{q}{p})$.

We turn  the hemispheres $\south$ and $\north$ into networks by equipping $\south$ with the  sources  $v(i)$, for $i\in V(T)$ and the sinks $e(\arcInput{p}{q})$, for $\arc{p}{q} \in A_0$, and $\north$ with the sources $v'(i)$, for $i\in V(T)$, and the sinks $e'(\arcInput{q}{p})$, for $\arc{p}{q} \in A_0$. 

Note that both $\south$ and $\north$ are instances of the construction $\HY$, not only as digraphs, but now also as networks.
Indeed, $\south$ is $\HY$ for $Y=T_A$ since
$\Ori{T_A}=A_0$, and $\north$ is $\HY$ for $Y=T_A'$ since $\Ori{(T_A')}=\{\arc{q}{p}\;|\; \arc{p}{q} \in A_0\}$.

\subsection{Non-intersecting paths in networks $\south$, $\north$ and $\route$}\label{sec: glueing iso} 
We glue pairs of non-intersecting path systems, one in each hemisphere, in order to obtain a non-intersecting path system on the route network. 

Let $\Lambda^S \in \NIP(\south)$ and $\Lambda^N \in \NIP(\north)$. 
For each $\gamma=\arc{p}{q}\in A_0$, let $\Lambda^S_\gamma$ be the path of $\Lambda^S$ with terminus $e(\arcInput{p}{q})$, and let $\Lambda^N_\gamma$ be the path of $\Lambda^N$ with terminus $e'(\arcInput{q}{p})$.  
Let $\Gamma(\Lambda^S, \Lambda^N)=\{\Lambda_\gamma^S \Psi(\Lambda_\gamma^N) \,|\, \gamma \in A_0 \}$. 

\begin{p}\label{prop: bij route = south x north}
The map $\Gamma$ is a bijection from  $\NIP(\south) \times \NIP(\north)$ to $\NIP(\route)$.
\end{p}
\begin{proof}
Clearly $\Gamma$ embeds  $\NIP(\south) \times \NIP(\north)$ into $\NIP(\route)$.

Let us check that this embedding is surjective.
Let $\Lambda \in \NIP(\route)$.
Any of the $n$ paths of $\Lambda$ has its origin in $\south$ and terminus in $\north$. Hence the path has a bridge as one of its steps. Since the path system is non-intersecting, each path passes through a different bridge.

Given $\gamma=\arc{p}{ q} \in A_0$, set 
 $\Lambda_{\gamma}$ to be the path of $\Lambda$ that passes through the bridge $(e(\arcInput{p}{q}), e'(\arcInput{p}{q}))$. Then $\Lambda_\gamma$ splits as $\mu \mu'$ where $\mu$ is a path in $\south$ with terminus $e(\arcInput{p}{q})$, and $\mu'$ is a path in $\north$ with origin $e'(\arcInput{p}{q})$.  Set $\Lambda_\gamma^S=\mu$ and $\Lambda_\gamma^N=\Psi(\mu')$. The pair $(\Lambda^S, \Lambda^N)$, where $\Lambda^S=\{\Lambda_\gamma^S\,|\, \gamma \in A_0\}$ and $\Lambda^N=\{\Lambda_\gamma^N\,|\, \gamma \in A_0\}$, has image $\Lambda$ under  $\Gamma$. 
\end{proof}

We claim that there is exactly one non-intersecting path in the route network. That is, that $\NIP(\route)$ has cardinality one. To prove this,  it will be enough to show  both $\NIP(\south)$ and $\NIP(\north)$ have cardinality one.

Since both hemisphere networks are instances of the construction $\NIP(\HY)$ (for $Y=T_A$ and $Y=T'_A$). Our task reduces to showing that   $\NIP(\HY)$  has only one element, for any left-right tree $Y$. 

\section{Contractions and liftings of paths}
\label{subsec: lift and proj}

Let $\treeU$ be a left-right tree. We showed in Proposition \ref{south is acyclic} that any walk  on $\hemi{Y}$ is a path.
In this section, we will show that any path in  $\hemi{Y}$ can be projected onto both the earring tree $\earring{Y}$ and the left-right tree $\treeU$. We will also show that paths in $\treeU$  can be lifted to paths in $\hemi{Y}$.

We define the \emph{contraction map} $\co$ on the vertices of $\earring{Y}$  that sends $M(i, j)$, $M_0(i)$, and $\M_0(i)$ to $i$, and $\co(r)=r$. 
Recall that the construction of $\earring{Y}$ involves a small parameter $\varepsilon > 0$, so that the embedding of every vertex $P$ in it depends on $\varepsilon$; we think of $\co(P)$ as $\lim_{\varepsilon \to 0} P(\varepsilon)$. 

We define the \emph{contraction $\co(\mu)$ of a walk} $\mu = P_1 P_2\cdots P_k$ on $\earring{Y}$ as the walk on $Y$ obtained by replacing in $\co(P_1) \co(P_2) \cdots \co(P_k)$ each one-letter maximal subword with a single occurrence of this letter (for instance, $11122234455$ is changed into $12345$).

\begin{ex}
Let $T_A$ be the tree of our running example. The unique path $W$ in $\hemi{T_A}$ from $v(4)$ to $e(2,6)$
corresponds (under the correspondence between non-trivial paths in $\HY$ and paths in $\earring{Y}$) to the path
\begin{multline*}
 w = \M_0{(4)}
M_0(4) M(4,7) M(4,9) M(4,1)
M(1,4) M_0(1)     \\   
M(1,r)
\,r\,
M(2,r) M(2,5) M_0(2) M(2,6) M(6,2).
\end{multline*}
See Figure \ref{fig: uniquepathv4e26} for an illustration of $w$.

\begin{figure}[h] \tiny \centering
\includegraphics[width = .4\textwidth, page = 28]{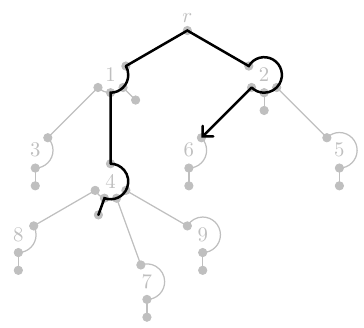}
 \caption{
 The path $w$ of  $\earring{T_A}$
 with first step $v(4)$ and last step $e(\arcInput{2}{6})$.}
  \label{fig: uniquepathv4e26}
\end{figure}
To compute the contraction of $ w $,  replace each occurrence of $M(i,j)$,  $M_0(i)$ and $\M_0(i)$ with $i$ obtaining $ 44444111r22226$. Thus, $\co(\mu)= 41r26$. 
\end{ex}

\begin{lem}[Lifting  of a path to $\HY$]\label{lem: lifting to hemi}
Let $\pi = i_0 i_1 \ldots i_k$ be a non-trivial path in $Y$.
Then, there exists a unique path in $\hemisphere$ from $v(i_0)$ to $e(\arcInput{i_{k-1}}{i_k})$. 
\end{lem}
 We call this path the \emph{lifting  to $\hemisphere$ of  $i_0i_1\cdots i_k$}.
 
\begin{proof}
Let $R$ and $Q$ be the vertices of $\earring{Y}$ such that $e(\arcInput{i_{k-1}}{i_k})=\arc{R}{ Q}$. 
Since the nodes of $\hemisphere$ are the arcs supported on $\earring{Y}$, 
it suffices to prove that there exists a path in $\earring{Y}$ whose first step is $v(i_0) = \arc{\M_0(i_0)}{M_0(i_0)}$ and whose last step is $e(\arcInput{i_{k-1}}{i_k}) = \arc{R}{Q}$.

Since $\earring{Y}$ is a tree, there exists a unique path $\mu=P_0 P_1 \ldots P_{\ell}$ in $\earring{Y}$ from $\M_0(i_0)$ to $Q$.
Thus $P_0=\M_0(i_0)$ and $P_{\ell}=Q$. 
Since $\M_0(i_0)$ has as unique neighbor $M_0(i_0)$, we have $P_1=M_0(i_0)$.

The contraction $\co(\mu)$ is a path in $Y$ from $i_0$ to $i_{k}$; it is therefore $\pi$.  
Let $m$ be the last index such that $\rho(P_m) = i_{k-1}$. Then $\rho(P_{m+1})=i_k$. 
Since $e(\arcInput{i_{k-1}}{i_k})$ is the only arc of $\earring{Y}$ from a vertex whose image under $\rho$ is $i_{k-1}$ to one whose image is $i_k$, we have $\arc{P_{m}}{P_{m+1}}=e(\arcInput{i_{k-1}}{i_k})=\arc{R}{ Q}$.
Since $\mu$ is a path (no repeated vertex) and $P_{m+1}$ and $P_{\ell}$ are both equal to $Q$, we have $m+1=\ell$ and thus $m=\ell-1$ and $P_{\ell-1}=P_m=Q$. 

The uniqueness follows from the unique walk property of $\HY$ established in Lemma \ref{south is acyclic}.
\end{proof}

\section{Catalysts induce path systems in  \texorpdfstring{$\route$}{R(A)}.} 
\label{subsec: Lambda lifts npaths}

 In this section, we show how every catalyst $(\sigma, f)$ in the arrowflow class $ C(A)$ induces a path system in the route network $\route$.

\begin{lem}[Lifting of a marked path to $\route$] \label{lifting marked path}
Let $(\pi;\gamma)$ be a marked path in $T_A$ such that $\gamma \in A_0$. Let $i$ and $j$ be the origin and terminus of $\pi$.
Then there exists a unique path in $\route$ from $v(i)$ to $v''(j)$ that admits $\arc{e(\gamma)}{ e'(\gamma})$ as a step.
\end{lem}
We call this path the \emph{lifting to $\route$ of $(\pi;\gamma)$.}  
By a slight abuse of notation, we also refer to the lifting of $\left(\pathTO{i}{j}; \gamma\right)$ (marked path of $T_A$) as the  lifting of $\left(\pathT{i}{j}; \theta^{-1}(\gamma)\right)$  (corresponding marked path of $T$).

 \begin{proof}
 Write $\pi$ as $i_0i_1\cdots i_m$, with $i_0=i$ and  $i_m=j$. There exists $k$ such that $\gamma=i_k i_{k+1}$.
 Let $ v(i_0) \cdots e(\arcInput{i_k}{i_{k+1}})$  be the lifting to $\south$ of the path $i_0 i_1 \cdots i_k i_{k+1}$. 
 Let $e'(\arcInput{i_k}{i_{k+1}}) \cdots v'(i_m)$ be the image under the  anti-isomorphism $\Psi$ of the lifting to $\north$ of the path $i_m i_{m-1} \cdots i_{k+1} i_k$.
 Then, since $\arc{e(\arcInput{i_k}{i_{k+1}})}{ e'(\arcInput{i_k}{i_{k+1}}})$ is a bridge of $\route$,  the concatenation
\[
v(i_0) \cdots e(\arcInput{i_k}{i_{k+1}}) e'(\arcInput{i_k}{i_{k+1}}) \cdots v'(i_m).
\]
is a path in $\route$ that 
fulfills the requirements of the lemma. The uniqueness follows from the unique walk property  of hemispheres $\south$ and $\north$ (Lemma \ref{south is acyclic}).
 \end{proof}
 
\begin{lem}[Lifting to $\route$ of a catalyst.]
\label{lem: lifting to route of a catalyst}
Given a catalyst $\kappa=(\sigma; f) \in C(A)$, there is a unique path system $\Lambda$ of $\route$ such that, for each $i$, the path of $\Lambda$ with origin $v(i)$ has terminus $v''(\sigma(i))$ and passes through the bridge $\arc{e(\rootins\circ f(i))}{e'(\rootins\circ f(i))}$.
\end{lem}

This path system is  the \emph{lifting to the route network $\route$} of catalyst $\kappa$. 
\begin{proof}
This is the path system whose paths are the liftings of the marked paths $\left(\pathT{i}{\sigma(i)}; f(i)\right)$.
The uniqueness follows from the uniqueness of the lifting of a marked path (\ref{lifting marked path}). 
\end{proof}

\begin{ex}\label{ex:lift}
Let  $\kappa=(\sigma, f)$ be the catalyst of Example \ref{example_signed_path}, defined by  the following table:
\[
\begin{array}{c|ccccccccc}
i & 1 & 2 & 3 & 4 & 5 & 6 & 7 & 8 & 9 \\\hline
\sigma(i) & 6 & 5 & 8 & 7 & 3 & 2 & 9 & 4 & 1 \\
f(i) & 12 & 25 & 31 & 47 & 21 & 62 & 49 & 84 & 41 
\end{array}
\]

Let $\route$ be the route network in Figure \ref{Fig: lift}. We illustrate the lifting  of $\kappa$ into a path system in $\route$ by lifting two of the paths explicitly.

The lifting $\Lambda_9$ to $\route$ of  $(\pathT{9}{1};\arc{4}{1})$ starts at $v(9)$, passes through the bridge $\arc{e(\arcInput{4}{1})}{ e'(\arcInput{4}{1})}$ and ends at $v''(1)$.  It corresponds to the path in $\earring{T_A}$ with first step $v(9)=\M_0(9)M_0(9)$ and last step $e(4,1)=M(4,1)M(1,4)$, followed by the path  in $\earring{T'_A}$ with first step $e'(4,1)=M'(4,1)M'(1,4)$ and last step $v''(1)=M'_0(1)\Mp_0(1)$.  See Figure \ref{fig: lift of one path}.

To obtain the  lifting $\Lambda_1$ of  $(\pathT{1}{6};12)$, we first lift $12$ to $r2\in A_0$. Then, $\Lambda_1$ starts at $v(1)$, passes through the bridge $\arc{e(\arcInput{r}{2})} {e'(\arcInput{r}{2})}$ and ends at $v''(6)$.
         
\begin{figure}[h]
    \centering\hfill
    \begin{subfigure}[b]{0.23\textwidth}
        \includegraphics[width = \textwidth, page=29]{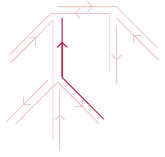}
        \caption{A marked path.}
    \end{subfigure}\hfill
    \begin{subfigure}[b]{0.76\textwidth}
        \includegraphics[width = \textwidth, page=30]{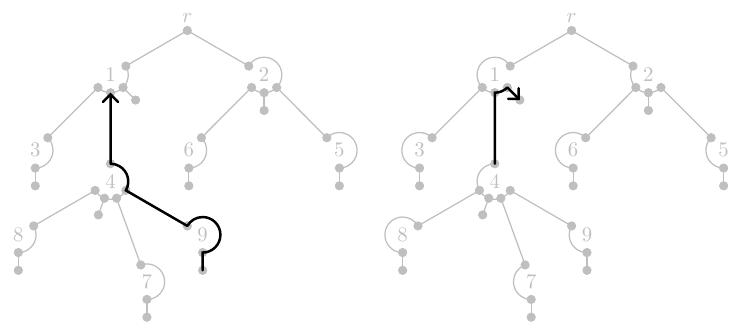}
        \caption{The paths in $\earring{T_A}$ and $ \earring{T_A'}$ corresponding to $\Lambda_9$.}
    \end{subfigure}\hfill
    \caption{Lifting $\Lambda_9$ of  
    $(\pathT{9}{1};41)$ to $\route$.}
    \label{fig: lift of one path}
\end{figure}

Figure \ref{Fig: lift} illustrates the path system $\{\Lambda_1, \ldots, \Lambda_n\}$ lifting catalyst $\kappa$. Paths $\Lambda_i$ are represented by the corresponding paths of $\earring{T_A}$ and $\earring{T'_A}$ through the node-arc correspondence of Lemma \ref{lem: proj from linegraph}.  
The direction of each path is implied, since $\Lambda_i$  goes from $\M_0(i)$ to some $\arc{P}{Q}$ in $\south$ and then from $\arc{Q'}{P'}$ to some $\Mp_0(j)$ in $\north$.

\begin{figure}[h]
    \centering\hfill
    \begin{subfigure}[b]{0.23\textwidth}
        \includegraphics[width = \textwidth, page=31]{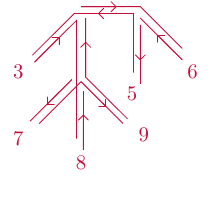}
        \caption{A catalyst $\kappa$.}
    \end{subfigure}\hfill
    \begin{subfigure}[b]{0.76\textwidth}
        \includegraphics[width = \textwidth, page=32]{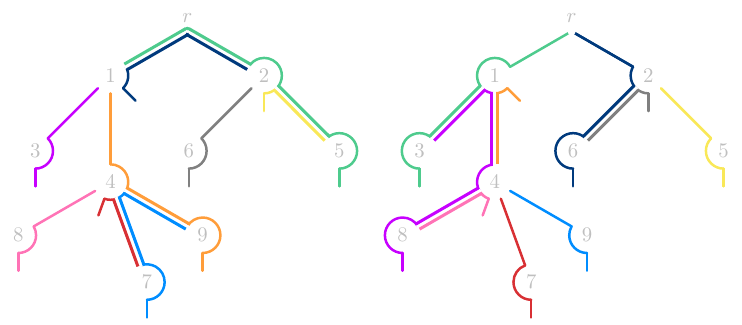}
        \caption{Its lifting $\Lambda(\kappa)= \{\Lambda_1, \ldots, \Lambda_n\}$.}
    \end{subfigure}\hfill
    \caption{Lifting of a catalyst to $\route$.}
    \label{Fig: lift}
\end{figure}
\end{ex}

\medskip

The lifting map for catalysts embeds $C(A)$ into $\PS(\route)$, the set of path systems of $\route$. However, it is not a surjection: path systems whose paths share some bridges can not be reached. 

We say that a  path system of $\route$ is {\em full} when every bridge $\arc{e(\gamma)}{e'(\gamma)}$ belongs to at least one (and thus: exactly one)  of its paths.

\begin{lem}\label{lem: lifting is bijective}
    The operation of lifting is a sign-preserving bijection from $C(A)$ to the set of full path systems in $\route$.
\end{lem}
\begin{proof}
Let $\Lambda$ be a full path system in $\route$. Let $\sigma$ be its underlying permutation. 
Define a map $\tilde{f} : V(T) \to A_0$ that sends vertex $i$ to the arc $\gamma$ corresponding to the bridge $\arc{e(\gamma)}{e'(\gamma)}$ of $\Lambda_i$. 
Then $(\sigma, \rootins^{-1} \circ \tilde{f})$  is a catalyst in $C(A)$ whose lifting is $\Lambda$.

Finally, since, by construction, the underlying permutation of a catalyst is the permutation induced by its lifting,
the lifting map is a sign-preserving  bijection.
\end{proof}
 
\begin{thm}\label{thm: full}
    Let $A$ be a unital arrowflow. Then,
    \begin{equation}\label{eq: CANIP}
 \sum_{\kappa\in C(A)} \sign(\kappa) =
\sum_{\Lambda \in \NIP(\route)} \sign(\Lambda).
    \end{equation}
\end{thm}
\begin{proof}
By Lemma \ref{lem: lifting is bijective},
\[
\sum_{\kappa\in C(A)} \sign(\kappa) =
\sum_{\Lambda\in\Full(\route)} \sign(\Lambda),
\]
where $\Full(\route)$ denotes the set of full path systems on $\route$.
Recall that $\route$ is acyclic by Corollary \ref{cor: paths in N and R}. Thus there is a  Lindström--Gessel--Viennot involution $\phi$ (Lemma \ref{lem: LGV}) on $\PS(\route)$. Since $\phi$ preserves the multiset  of steps of the path systems, it sends full path systems to full path systems. Hence, $\phi$ restricts to a sign-reversing involution of $\Full(\route)$, whose fixed points are the non-intersecting path systems that are full. But all non-intersecting path systems are full, whence \eqref{eq: CANIP}.  
\end{proof}
\section{The unique non-intersecting path system of \texorpdfstring{$\hemi{\treeU}$}{H(Y)}}
 \label{NonintersectingHemisphere}

Let  $Y$ be a left-right tree, whose root $r$ has two children, either both left or both right children. Under this assumption, the results of this section can be applied both  to $T_A$, where both children of the root are right children, and to $T'_A$, where they are both left children.

In Section~\ref{technical proofs}, we show that the hemisphere network $\hemi{Y}$ has a unique non-intersecting path system. Assuming this result, we describe in this section the unique path system  using the classical Depth-First Search (DFS) walk on a plane-rooted tree and the inorder traversal of its vertices.

Given a plane rooted tree $U$ with root $r$,  its {\em {Depth-First-Search} (DFS) walk} is the walk $\dfs(r)$, where $\dfs(i)$ is the walk defined recursively for each vertex $i$ as follows:
if $i$ is a leaf then $\dfs(i)$ is the one-term sequence $(i)$; 
otherwise, let $j_1 <_i j_2 <_i \cdots <_i j_k$ be the children of $i$, then
\[
\dfs(i) = 
i 
\;  \dfs(j_1)
\; i
\;\dfs(j_2)
\; i
\cdots 
i
\; \dfs(j_k)
\; i.
\]
The DFS walk of $U$ is a closed walk with length  $2m$, where $m$ is the number  of vertices of $U$.
When convenient, we also consider the DFS walk of $U$ as indexed by the integers modulo $2m$;
in that case we refer to it as the \emph{cyclic DFS walk of $U$}.

We will make use of the following property of the DFS walk and the cyclic DFS walk:
 every arc supported on $U$ appears exactly once as a  step of it (\cite[Lemma 3.2]{Even}).

\begin{ex}\label{ex: DFS walk}
    The DFS walk for the plane rooted tree $T_A$ of Figure  \ref{fig: DFS on T0} is $r13148474941r26252r$. 
\begin{figure}[h!]
\begin{subfigure}[b]{0.45\textwidth}
\begin{center}
        \includegraphics[width = 6cm, page=33]{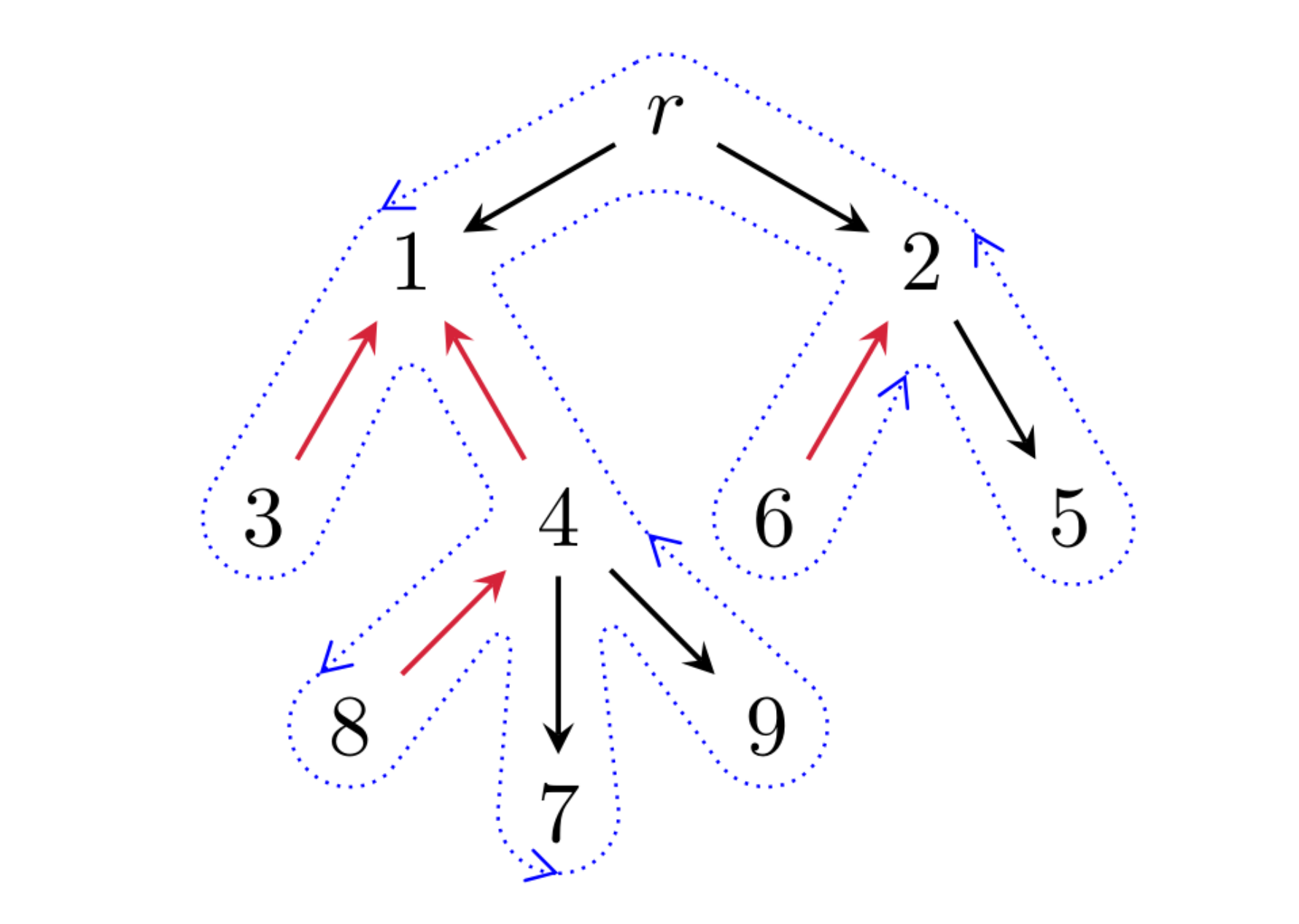}
\end{center}
\caption{The cyclic DFS walk on the left-right tree $Y$ follows the dotted loop.\phantom{xxxxx}\phantom{xxxxx}\phantom{xxxxx}\phantom{xxxxx}\phantom{xxxxx}\phantom{xxxxx}
}\label{fig: DFS on T0}
\end{subfigure}
\hfill
\begin{subfigure}[b]{0.44\textwidth}\centering
        \includegraphics[width = \textwidth, page=34]{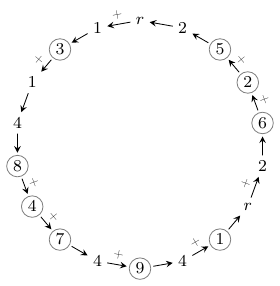}
\caption{The same cyclic DFS walk, unfolded, with intermediate visits circled, and steps in $\OriY$ marked  ``+''.}
\label{fig: W}
\end{subfigure}
\caption{}
\label{fig: DFS walk}
\end{figure}
\end{ex}
The classical {\em inorder} of the vertices of a binary tree (see for instance \cite[Ch.~4]{Sedgewick}) extends straightforwardly to our left-right trees. 
For each vertex $j$ of a left-right tree $Y$, call the \emph{intermediate visit} to $j$ its visit ($k$ such that $i_k=j$) that takes place after all visits to all its left children in the DFS walk $i_0 i_1 \cdots i_{2m}$ of $Y$, but before all visits to all its right children. The inorder on the vertices of $Y$ is the order of their intermediate visits.

\begin{ex}\label{ex: inorder}
  For the left-right  tree of  Figure \ref{fig: DFS on T0}, the intermediate visits are the circled entries in
  $  \circled{r}1\circled{3}14\circled{8}\circled{4}\circled{7}4\circled{9}4\circled{1} r2\circled{6}\circled{2}\circled{5}2 r$.
  Thus, the  vertices in inorder are $r384791625$. 
\end{ex}

In the following theorem, we say that an arc $\arc{i}{j}$ of $Y$ is \emph{upwards} if $j$ is the parent of $i$, and \emph{downwards} if $j$ is a child of $i$.

\begin{thm}\label{description of the NI in HY}
Let $Y$ be a left-right tree, whose root has either two right children or two left children.
Then the hemisphere network $\hemi{Y}$ has exactly one non-intersecting path $\LU$.

Let $i_1, i_2, \ldots, i_n$ be the elements of $\pV(Y)$ in inorder, and let $\sigma_Y$ be the cyclic permutation $(i_1 i_2 \cdots i_n)$.  For each $i \in \pV(Y)$, let $\mu_i$ be the path from $i$ to $\sigma_Y(i)$. 
Then:
\begin{myenum}
\item  \label{HY item1}
The paths $\mu_i$ are subwalks of the cyclic DFS walk on $Y$, and each arc supported on $Y$ is a step of exactly one of them. 
\item \label{HY item2} Each path $\mu_i$ has exactly one step belonging to $\OriY$.  Call it $h_Y(i).$
\end{myenum}
Fix $i\in \pV(Y)$. Write $\mu_i=j_0 j_1 \cdots j_k j_{k+1} \cdots j_m$, with $j_k j_{k+1}=h_Y(i)$ (the step belonging to $\OriY$). Then:
 \begin{myenum}
 \setcounter{enumi}{2}
 \item \label{HY item3} The path of $\LU$ with origin $v(i)$  is the lifting to $\HY$ of $j_0 \cdots j_k j_{k+1}$.
\item \label{HY item4} All steps in $j_0 \cdots j_k$ are upwards, and all steps in $j_{k+1} \cdots j_m$ are downwards.  
\end{myenum}
\end{thm}

\begin{ex}
  In Example \ref{ex: inorder}, the vertices different from the root in inorder are
    $384791625$. 
    The permutation $\sigma_{T_A}$ is the cycle $(384791625)$. 
    The paths $\mu_i=\pathTO{i}{\sigma_{T_A}(i)}$ are $3148$, $84$, $47$, $749$, $941$, $1r26$, $62$, $25$ and $52r13$. 
    The paths of $\LU$ are the liftings of $31$, $84$, $47$, $749$, $941$, $1r2$, $62$, $25$ and $52r1$.
    \end{ex}

\section{The unique non-intersecting path system of \texorpdfstring{$\route$}{R(A)}}
\label{NonintersectingRoute}

This section wraps up our combinatorial argument. We show that there is a unique non-intersecting path system in $\route$, and that its underlying permutation is always an $n$-cycle. 
This amounts to the remaining half of Theorem \ref{thm: sum on class}. 
In this section, we let $A$ be a fixed unital arrowflow.

\begin{lem}\label{perm inverse}
The permutation $\sigma_{T_A'}$, associated to the mirror image $T_A'$ of $T_A$, is the inverse of $\sigma_{T_A}$.
\end{lem}
\begin{proof}
Let $i_1, i_2 ,\ldots, i_n$ be the vertices of $T$ in inorder with respect to $T_A$. The DFS walk on $T'_A$ is the DFS walk on $T_A$ reversed. On the other hand, the left  (resp. right) children in $T'_A$ are the right (resp. left) children in $T_A$. It follows that the vertices of $T$, in inorder with respect to $T_A'$,  are   $i_n,\ldots, i_2, i_1$. 
 Then $\sigma_{T_A}=(i_1 i_2 \cdots i_n)$ and  $\sigma_{T_A'}=(i_n \cdots i_2 i_1)$. These two cycles are inverse of each other. 
\end{proof}

Remember that after Theorem  \ref{description of the NI in HY}\eqref{HY item2}, each path $\pathTO{i}{\sigma_{T_A}(i)}$ has a unique step $h_{T_A}(i)$ in $A_0$.

\begin{p}\label{prop: exists unique n-path}
There is exactly one non-intersecting path system in the route network $\route$. 
Its paths are the liftings to $\route$ of the marked paths 
\[
\left(\pathTO{i}{\sigma_{T_A}(i)}; h_{T_A}(i)\right)
\]
of $T_A$, for $i \in V(T)$. 
Its underlying permutation is the cycle $\sigma_{T_A}$. 
\end{p}

\begin{proof}
The existence and uniqueness follow straightforwardly from the existence of a bijection $\NIP(\route) \cong \NIP(\south) \times \NIP(\north)$ (Proposition \ref{prop: bij route = south x north}) and the fact that each of $\NIP(\south)$ and $\NIP(\north)$ has exactly  one element (Lemma \ref{lem: unique in HY}).

Index the paths of the unique path systems $\Lambda$, $\Lambda^S$ and $\Lambda^N$ of $\route$, $\south$ and $\north$ with the $\gamma \in A_0$ as in Section \ref{sec: glueing iso}.

Let $\gamma\in A_0$. 
By Theorem \ref{description of the NI in HY}\eqref{HY item1} and \eqref{HY item2}, applied to $Y=T_A$, there exists a unique $i\in V(T)$ such that $\gamma$ is a step of $\pi=\pathTO{i}{\sigma_{T_A}(i)}$, and $\gamma=h_{T_A}(i)$. 
Write $\pi$ as $j_0\cdots j_k j_{k+1} \cdots j_m$, with $j_k j_{k+1}=\gamma$. 
Then $\Lambda^S_{\gamma}$ is the lifting to  $\south$ of $j_0\cdots j_k j_{k+1}$ by Theorem \ref{description of the NI in HY}\eqref{HY item4}. In particular its origin is $v(i)$.

By Lemma \ref{perm inverse}, $i=\sigma_{T_A'}(j)$ for $j=\sigma_{T_A}(i)$. 
Note also that $\Psi(\gamma)$ is a step of $\Psi(\pi)=\pathTO{\sigma_{T_A}(i)}{i}$, which is $\smallpathTO{j}{\sigma_{T_A'}(j)}$.
Applying Theorem \ref{description of the NI in HY} to $Y=T'_A$, we get that  $\Lambda_{\gamma}^N$ is the lifting of $j_m \cdots j_{k+1} j_k$ to $\north$.  
In particular, the origin  of $\Lambda_{\gamma}^N$ is $v'(j)$.  
Thus, the terminus of $\Psi(\Lambda^N_{\gamma})$ is $v''(j)$.

The path $\Lambda_\gamma=\Lambda_{\gamma}^S \Psi(\Lambda^N_{\gamma})$  has origin $v(i)$, terminus $v''(j)$ and has the bridge $e(\gamma)e'(\gamma)$ as a step. 
Therefore, by Theorem \ref{lifting marked path}, $\Lambda_{\gamma}$ is the lifting to $\route$ of the marked path $(\pathTO{i}{j}; \gamma)$. 

This holds for all $n$ elements $\gamma$ of $A_0$. 
This proves that the paths of $\Lambda$ are the liftings of the $n$ marked paths $\left(\pathT{i}{\sigma_{T_A}(i)}; h_{T_A}(i)\right)$. 
\end{proof}

This shows the second half of Theorem \ref{thm: sum on class}.

\begin{cor}\label{thm: unital sum}
      \[
    \sum_{\kappa \in C(A)} \sign(\kappa) = (-1)^{n-1}.
    \]
\end{cor}

\begin{proof}
    By Theorem \ref{thm: full}, 
    $
     \sum_{\kappa\in C(A)} \sign(\kappa)= \sum_{\Lambda \in \NIP(\route)} \sign(\Lambda).
   $
    By Proposition \ref{prop: exists unique n-path}, the latter sum has a unique summand, which is $(-1)^{n-1}$ since its permutation $\sigma_{T_A}$ is a cycle of length $n$.
\end{proof}

Through the lifting maps, the LGV involution on the path systems of $\route$
induces an involution $\varphi_A$ on the unital arrowflow class $C(A)$. By abuse of language, we still refer to $\varphi_A$ as the LGV involution on $C(A)$. 

We proceed to describe its unique fixed point.
Partition set $\arcs{T}\setminus A$, into two parts: those oriented upwards and those oriented downwards.  More precisely, let 
 \begin{align*}
 \Aup &= \{\arc{i}{j} \in \arcs{T} \; |\; \arc{i}{j} \not\in A \text{ and $j$ is the parent of $i$ in $T_A$} \},\\
 \Adown &= \{\arc{i}{j} \in \arcs{T} \; |\; \arc{i}{j} \not\in A \text{ and $j$ is a child of $i$ in $T_A$} \}.
 \end{align*}
Thus 
\(\arcs{T} = A \sqcup \Aup \sqcup \Adown \), where $\sqcup$ denotes the disjoint union. 

\begin{thm} \label{surviving catalyst}
For each vertex $i$ of $T$, let $\pi_i = \pathT{i}{\sigma_{T_A}(i)}$.
Then:
\begin{myenum}
\item\label{kA i} Each arc supported on $T$ is a step of exactly one of the paths $\pi_i$. 
\item\label{kA ii} Each of the paths $\pi_i$ has a unique step in $A$. Call it $f_A(i)$. 
\item\label{kA iii} The catalyst $\kappa_A=(\sigma_{T_A}, f_A)$ is the unique fixed point of the LGV involution on $C(A)$. 
\item\label{kA iv} For each $i$, all steps in the tail of the marked path $(\pi_i; f_A(i))$ are in $U(A)$, and all steps in its head are in $D(A)$. \end{myenum}
\end{thm}

\begin{ex}
  For the left-right tree of Figure \ref{fig: plane rooted tree C},  the fixed point of $\varphi_A$ is the catalyst $(\sigma_{T_A}, f_A )$ corresponding to the family of marked paths $\underline{31}48$, $\underline{84}$, $\underline{47}$, $7\underline{49}$, $9\underline{41}$, $\underline{12}6$, $\underline{62}$, $\underline{25}$ and $5\underline{21}3$.  
  Its permutation is $\sigma_{T_A}=(384791625)$, and $f_A$ maps each $i$ to the marked step of the marked path starting at $i$: $3$ to $\arc{3}{1}$, $8$ to $\arc{8}{4}$, $4$ to $\arc{4}{7}$, $7$ to $\arc{4}{9}$, and so on.
\end{ex}  

\begin{proof}
 Each of the $n$ paths  $\pathTO{i}{\sigma_{T_A}(i)}$ is equal to $\pathT{i}{\sigma_{T_A}(i)}$, except for the two of them that have $ra$ and $rb$ as a step.

 Let us consider the case of the path $\mu_i=\pathTO{i}{\sigma_{T_A}(i)}$ that has a step $ra$. Then the predecessor of $r$ in this path is $b$, and $\mu_i$ factorizes as $w bra w'$, with $ra$ its unique step in $A_0$. 
 Then  $\pathT{i}{\sigma_{T_A}(i)}$ factorizes as $w ba w' $, with $ba$ its unique step in $A$. The case of the path that has $rb$ as a step is similar.

 After this observation, the theorem follows by applying 
Theorem \ref{description of the NI in HY} and Proposition \ref{prop: exists unique n-path}.
\end{proof}

\section{Deformations of the Graham--Pollak formula}\label{sec: generalizations} 

The combinatorial framework developed to establish a combinatorial proof for the Graham-Pollak formula also extends to its generalizations. This framework further yields a new generalization, stated in Theorem \ref{main: formula indep}, from which all such results can be derived.

\subsection{The marked distance matrix}
Recall the definition of the marked distance $d'(i,j)$ from the introduction: we began by associating to each marked path 
$(\pi; \gamma)=(i_0~i_1\ldots i_{m}; i_{k} i_{k+1})$
in $T$ the monomial
\[
w(\pi;\gamma) 
=x_{i_0i_1}\cdots x_{i_{k-1}i_{k}}y_{i_k i_{k+1}}z_{i_{k+1}i_{k+2}}\cdots z_{i_{m-1}i_{m}}\]
in three families of commuting variables $x_{\gamma}$, $y_{\gamma}$, $z_{\gamma}$ attached to the arcs $\gamma$ supported on $T$.

Given two vertices $i$ and $j$, and $\pathT{i}{j}=i_0 i_1 \cdots i_d$ the path from $i=i_0$ to $j=i_d$, 
the {\em marked distance} between $i$ and $j$ was defined as 
\[
d'(i,j) = \sum_{p=0}^{d-1} w\!\left(i_0i_1 \cdots i_d; \arc{i_{k}}{i_{k+1}}\right).
\]
Finally, the {\em marked distance matrix} $M'(T)$ was defined as $(d'(i,j))_{1\le i,j\le n}$.

The \emph{weight of a set of marked paths} is defined as the product of the weights of its paths. 
This applies to catalysts, after identifying them with their corresponding set of marked paths: the \emph{weight of a catalyst} $(\sigma, f)$ is thus defined as 
\[
w(\sigma, f) = \prod_{i\in V(T)} w\!\left(\pathT{i}{\sigma(i)}; f(i)\right).
\]
We define the \emph{tail} (\emph{head}) \emph{of a catalyst} to be the multiset of tails (heads) of its marked paths.

Our results generalize straightforwardly  to the weighted setting whenever the sign-reversing involutions defined on the arrowflow classes $C(A)$ are all weight-preserving.

\subsection{Zero-sum involutions}

Let us show that  the involutions $\varphi_A$ on classes of  zero-sum catalysts are weight-preserving, whenever a simple condition on the weights ($x_{ij}x_{ji}=1$) is fulfilled.

\begin{lem}
    Let $A$ be a disconnected zero-sum arrowflow. Assume $x_{ij} = x_{ji}^{-1}$ for all $\arc{i}{j}\in \arcs{T}$.
    The involution $\varphi_A$ defined in Lemma \ref{prop zero-sum disconnected} is weight-preserving.
\end{lem}
\begin{proof}
    Let $\kappa = (\sigma, f)\in C(A)$, and let $i$ and $j$ as in Lemma \ref{prop zero-sum disconnected}.
    
    If  $\{i,j\}$ does not separate $\sigma(i)$ and $\sigma(j)$ (as in Figure \ref{fig: disconnected 1}), 
    $\kappa$ and $\varphi_A(\kappa)$ have the same head and tail, and thus the same weight.
    
    If   $\{i,j\}$ separates $\sigma(i)$ and $\sigma(j)$ (as in Figure \ref{fig: disconnected 2}), then  $\varphi_A(\kappa)$ and $\kappa$ have the same  head. The 
    tail of $\varphi_A(\kappa)$ is obtained from the tail of $\kappa$ by either removing or adding one occurrence of each of $\arc{i}{j}$ and $\arc{j}{i}$. Accordingly, the weight of $\varphi_A(\kappa)$ is obtained from that of $\kappa$ by multiplying or dividing by $x_{ij} x_{ji}$. This has no effect since $x_{ij}x_{ji}=1$.
\end{proof}

\begin{lem}
Let $A$ be a connected zero-sum arrowflow.
    The  involution $\varphi_A$ defined in Lemma \ref{prop zero-sum connected} is weight-preserving.
\end{lem}
\begin{proof}
    The involution $\varphi_A$ just swaps the tails of two of the marked paths of the catalyst (see Figure \ref{fig: double edge}). 
    In total, it preserves the head  and tail of the catalyst. Hence, it is weight-preserving.
\end{proof}

As a result, the determinant of the marked distance matrix is the signed, weighted sum of all unital catalysts,
\[
\det M'(T) = 
\sum_{A~\text{unital}}~
\sum_{\kappa\in C(A)}\sign(\kappa)w(\kappa).
\]

\subsection{Lifting to the weighted Route Networks}

Recall that the sign-reversing involution $\varphi_A$ on $C(A)$ was defined as the
involution induced on $C(A)$ by the
LGV involution on the path systems of $\route$ (through the lifting maps).
In order to prove that $\varphi_A$ is weight-preserving, we define weights for the path systems of $\route$, such that the LGV involution, as well as the lifting of catalysts, are weight-preserving.

For each unital arrowflow $A$, consider the plane rooted tree $T_A$, and its orientation $A_0$ obtained by substituting $\arc{a}{b}$ and $\arc{b}{a}$ for $\arc{r}{a}$ and $\arc{r}{b}$ (as in Subsection \ref{ConstructionT0}). 

The lifting of catalysts to $\route$ is based on the lifting of marked paths of $T$ to paths of $\route$, which is done  in two steps: firstly, marked paths of $T$ are lifted to marked paths of $T_A$; second, marked paths of $T_A$  are lifted to paths of $\route$.
We will show that we can equip the marked paths of $T_A$ and the paths of $\route$ with weights so that the liftings preserves the weights. 

For the lifting to $T_A$: we extend the weight defined on marked paths of $T$ to marked paths of $T_A$ by assigning weights to the new arcs as follows: to $\arc{r}{a}$ the same weights as $\arc{b}{a}$ ($x_{ra}=x_{ba}$ and likewise for $y$ and $z$), to $\arc{r}{b}$ the same weights as $\arc{a}{b}$ ($x_{rb}=x_{ab}$ \ldots) and to $\arc{a}{r}$ and $\arc{b}{r}$ all weights $1$ ($x_{ar}=y_{ar}= \cdots = z_{br}=1$).

\begin{lem}
  The lifting of marked paths from $T$ to $T_A$ (Section \ref{sec: root insertion}) preserves the weights.
\end{lem}  

\begin{proof}
  A marked path of $T$ is different from its lifting to $T_A$ only if has a step $ab$ or  a step $ba$.

  Consider a marked path $(\pi; \gamma)$ of $T$ having $\arc{a}{b}$ as a step. Let $(\pi_0; \gamma_0 )$ be its lifting to $T_A$.
  The path $\pi_0$  is obtained from $\pi$ by replacing $ab$ with $arb$.
  If this replacement takes place in the tail, then the weight is changed by replacing a factor $x_{ab}$ with $x_{ar} x_{rb}$.
  This has no effect since $x_{ar}=1$ and $x_{rb}=x_{ab}$.
  Similarly, if the replacement takes place in the head.
  Finally, if $\arc{a}{b}$ is the marked step $\gamma$,
  then $\pi$ factorizes as $\mu a b \mu'$, with $\mu a$ its tail and $b \mu'$ its head, while $\pi_0=\mu a r b \mu'$ with $\mu a r$ its tail and $b\mu'$ its head. The new marked step is $\arc{r}{b}$, instead of $\arc{a}{b}$.
  In the weight, a factor   $y_{ab}$ is replaced  with $y_{rb}$ and a new factor $x_{ar}$ appears.   
   The product is unaffected since $x_{ar}=1$ and $y_{rb}=y_{ab}$.

    The case of a step $\arc{b}{a}$ is treated similarly. 
\end{proof}

We equip $\route$ with a weight  function $w$ on its arcs by assigning
\begin{myenum}
\item  to all arcs of $\mathcal{S}$ leaving  some node $e(\gamma)$, the weight $x_{\gamma}$;
\item   to each bridge $\arc{e(\gamma)}{ e'(\gamma})$, the weight $y_{\gamma}$;
    \item to all arcs of $\mathcal{N}$ arriving to some node  $e'(\gamma)$, the weight $z_{\gamma}$;
    \item to any other arc, the weight  $1$.
\end{myenum}
The weight of a path in $\route$ is defined as the product of the weights of its steps.

\begin{lem} 
  The lifting of marked paths of $T_A$ with marked step in $A_0$, to paths of $\route$, defined in Lemma \ref{lem: lifting to route of a catalyst}, preserves the weights.
 \end{lem}

\begin{proof}
  Consider a marked path $j_0 j_1 \cdots j_k j_{k+1} \cdots j_m$ of $T_A$, with $j_k j_{k+1}$ its marked step. The nodes of the form $e(\gamma)$ and $e'(\gamma)$ in its lifting are, in order:
    $e(j_0j_1 ), e(j_1 j_2),\ldots , e(j_{k-1} j_k)$, 
    then the bridge edge $e(j_k j_{k+1}) e'(j_k j_{k+1})$, and then $e'(j_{k+1} j_{k+2}) ,\ldots,  e'(j_{m-1} j_m)$.    
Therefore the weight of the lifting is  
\[
x_{j_0j_1}\cdots x_{j_{k-1}j_{k}}y_{j_k j_{k+1}}z_{j_{k+1}ji_{k+2}}\cdots z_{j_{m-1}j_{m}}
\]
which is indeed the same as the weight of the marked path.
\end{proof}  

Define  the \emph{weight of a path system of $\route$} as the product of the weights of its $n$ paths. 
Since the weight  of a catalyst is   the product of the weights of the corresponding $n$ marked paths as well, and the lifting consists in lifting each marked path, the lifting of catalysts from $T$ to $\route$ is weight-preserving as well.

The LGV involution on the path systems of $\route$ is weight-preserving since it preserves the multiset of steps of the $n$ paths.

We conclude that the sign-preserving involution $\varphi_A$ on the unital class $C(A)$ is weight-preserving.

\subsection{The formula}

In each unital class $C(A)$, there is a unique catalyst $\kappa_A$, with sign $(-1)^{n-1}$, that is a fixed point of the LGV involution. Therefore
\[
\det M'(T) = \sum_{A \text{ unital}} (-1)^{n-1} w(\kappa_A).
\]
It follows from Theorem \ref{surviving catalyst} that 
\[
w(\kappa_A) = 
 \prod_{\gamma\in U(A)}
 x_\gamma
\quad \cdot 
 \prod_{\delta\in A}
 y_\delta
\quad  \cdot 
\prod_{\eta\in D(A)}
z_\eta
,
\]
where the sets $U(A)$,  and $D(A)$ denote the sets of arcs of $T$ that are upwards and downwards in $T_A$ and do not belong to $A$, as in Theorem \ref{surviving catalyst}.
The  sign of $\kappa_A$ is $(-1)^{n-1}$ by Proposition \ref{prop: exists unique n-path}.

For each edge $e = \{a,b\}$ of $T$, define $U(e)$ to be the set of arcs supported on $T$, different from  $\arc{a}{b}$ and $\arc{b}{a}$, that are upwards in the rooted tree obtained from $T$ by inserting a root $r$ by subdivision of $e$. 
Under the hypothesis $x_{ij}=x_{ji}^{-1}$ for all arcs $\arc{i}{j}$, the determinant of the marked distance matrix admits the compact expression presented in Theorem 
\ref{main: formula Emmanuel}, that we recall here and  will now derive:
\begin{equation}\label{eq: Theorem B recalled}
(-1)^{n-1} \sum_{e = \{a, b\}\in E} y_{ab}y_{ba} \prod_{\arc{i}{j} \in U(e)} (y_{ij} x_{ji} + y_{ji} z_{ij}).
\end{equation}

\begin{proof}[Proof of Theorem \ref{main: formula Emmanuel}]
   We have
  \[
  \det M'(T) = 
  \sum_{A \text{ unital}}   (-1)^{n-1} w(\kappa_A)
  \]
  with
\[    
w(\kappa_A) = 
  \prod_{\gamma\in U(A)}  x_\gamma
  \quad \cdot 
 \prod_{\delta\in A}  y_\delta
  \quad \cdot 
 \prod_{\eta\in D(A)} z_\eta 
.
\]

Let $A$ be a unital arrowflow and $e=\{a, b\}$ its marked edge. Then 
 \[
U(A) \sqcup A \sqcup D(A)
=
\arcs{T}
=  \{\arc{a}{b}, \arc{b}{a}\} \sqcup 
\bigsqcup_{\arc{i}{j} \in U(e)} \{\arc{i}{j}, \arc{j}{i}\},
 \]
For any  $\arc{i}{j}$ in $U(e)$, either it is in $A$, and then $\arc{j}{i} \in D(A)$, or it is in $U(A)$, and then $\arc{j}{i}\in A$. 
  We get therefore
 \[
 w(\kappa_A)  
   =  y_{ab} y_{ba}  \prod_{\arc{i}{j} \in U(e) \cap A} 
   y_{ij} z_{ji}
   \prod_{\arc{i}{j} \in U(e)\setminus A}
   x_{ij} y_{ji}
\]

Consider now all unital arrowflows with a given marked edge $e=\{a,b\}$. 
Any such arrowflow $A$ is determined by the set $B=U(e) \cap A$ of its upwards arcs, which can be any subset of $U(e)$: the other arcs in $A$ are the $\arc{j}{i}$ such that $\arc{i}{j} \in U(e)\setminus A$, and $\arc{a}{b}$ and $\arc{b}{a}$.
We get
\[
\det M'(T) =
(-1)^{n-1} 
\sum_{\{a,b\}\in E} y_{ab}y_{ba}
\sum_{B\subset U(e)} 
\prod_{\arc{i}{j} \in B} 
   y_{ij} z_{ji}
   \prod_{\arc{i}{j} \in U(e)\setminus B}
   x_{ij} y_{ji}
\]
which is indeed the expansion of  \eqref{eq: Theorem B recalled}. 
\end{proof}

For each edge $e$ of $T$, let $e^+$ and $e^-$ be the two arcs it supports.
Similarly, for any arc $\gamma$ of $T$, let $\gamma^-$ be its reverse.
Formula \eqref{eq: Theorem B recalled} can be rewritten as
    \begin{equation}\label{eq: Theorem B rewritten}
    \det M'(T) = (-1)^{n-1} \sum_{e\in E} y_{e^+}y_{e^-} \prod_{\gamma \in U(e)} (y_{\gamma} x_{\gamma^-} + y_{\gamma^-} z_{\gamma}).
    \end{equation}
    
\subsection{Generalization independent of tree structure}
The above generalization of the distance makes the determinant depend on the tree structure, and, in this sense, is not a generalization of the result of Graham and Pollak since it missed the important property that 
    ``the determinant of the distance matrix is independent of the tree structure.''

A simple way to further specialize the formula to make it independent of the tree structure is to impose
\[
(y_{\gamma} x_{\gamma^-} + y_{\gamma^-} z_{\gamma})
=
(y_{\gamma^-} x_{\gamma} + y_{\gamma} z_{\gamma^-}),
\]
so that $\gamma \in U(e)$ contributes to the formula as much as $\gamma^- \!\not\in U(e)$, for all $\gamma$.

Then we can introduce,  for each $e\in E$, a new variable $\alpha_{e}$ fulfilling
\[\alpha_{e}=\frac{y_{e^+}}{z_{e^+}-x_{e^+}} = \frac{y_{e^-}}{z_{e^-}-x_{e^-}}.
\]
Equivalently, $y_{e^+} = \alpha_{e} (z_{e^+}-x_{e^+})$ and $y_{e^-} = \alpha_{e} (z_{e^-}-x_{e^-})$.

\begin{proof}[Proof of Theorem \ref{main: formula indep}]
    It is clear that the generalized distance matrix $M_G(T)$ defined in the introduction is obtained from $M'(T)$ by specializing $y_{e^+}$ to $\alpha_{e} (z_{e^+}-x_{e^+})$ and $y_{e^-}$ to $\alpha_{e} (z_{e^-}-x_{e^-})$, for each edge $e$ of $T$.
    Specialize \eqref{eq: Theorem B rewritten} accordingly 
    to obtain for $\det M_G(T)$ the expression 
    \[
    (-1)^{n-1} \sum_{e\in E}
    \alpha_{e}^2 (z_{e^+}-x_{e^+})(z_{e^-}-x_{e^-})
    \prod_{\substack{f\in E\\f\neq e}}
    \alpha_{f}(z_{f^+}z_{f^-} - 1).
     \]
     Set now $x_e$ for $x_{e^+}$. Then $x_{e^-}=x_e^{-1}$.
\end{proof}

In this same spirit, we specialize $x_{e} = 1$ to obtain \cite[Thm.~A, case $x=0$]{CK19}.
This, in turn, implies every other additive generalization of the Graham--Pollak formula found in the literature (see \cite{CK19} for details).
We reiterate that ours is a fully combinatorial proof, in contrast to the more algebraic proofs found in the original paper.

\begin{cor}[Choudhury--Khare \cite{CK19}]\label{thm: CK}
Associate to each arc $\gamma$ supported on $T$ a variable $z_{\gamma}$, and to each edge $e$ of $T$ a variable $\alpha_e$. 

    Define the CK-weight of a marked path $(\pi;\gamma)$ to be
    \[
    \alpha_{e}(z_\gamma-1) \prod_{\delta\in \head({\pi;\gamma})} z_\delta.
    \]
    where $e$ is the edge that supports $\gamma$, and $\head({\pi;\gamma})$ is the set of steps of the head of $(\pi;\gamma)$.
        Define the CK-distance from $i$ to $j\in V(T)$ to be the  sum of the CK-weights of all marked paths from $i$ to $j$. 
        
    The determinant of the CK-distance matrix of $T$  is
    \[
    \left(\prod_{f\in E} \alpha_f (1-z_{f^+} z_{f^-})\right) \sum_{e\in E} \frac{\alpha_e (z_{e^+}-1)(z_{e^-} - 1)}{1-z_{e^+} z_{e^-}}.
    \]
\end{cor}
\begin{proof}
    Note that the CK-weight of a marked path is obtained from the weight of a marked path of Equation \eqref{eq: simple metric} by letting $x_{e} = 1$ for all $e$. 
    Theorem \ref{main: formula indep} then gives a formula for the determinant of the CK-distance matrix as
    \[
    (-1)^{n-1}
    \sum_{e\in E} \alpha_e^2
    (z_{e^+}-1)(z_{e^-}-1)
    \left(\prod_{f\ne e}\alpha_f(z_{f^-}z_{f^+}-1)\right). \qedhere
    \]
\end{proof}

\subsection{A remark on formal $q$-integers}

Among the deformations of the Graham--Pollak Formula, the $q$-analogue with weights $u_{\gamma}$ on  arcs \cite{LSZ} deforms 
the distance $d(i,j)$ into
\[
d_q(i,j) =  [u_{i_0 i_1}]
+ 
[u_{i_1 i_2}]
+ 
\cdots 
+
[u_{i_{d-1} i_d}]
\]
where $i_0 i_1 \cdots i_d$ is the path from $i$ to $j$, and $[u_{\gamma}]$ stands for $(q^{u_{\gamma}}-1)/(q-1)$.
The determinant of the matrix of the $d_q(i,j)$ is shown to be equal to \cite[Thm.~3]{LSZ}:
    \[
    (-1)^{n-1} 
    \sum_{e\in E} 
    [u_{e^+}][u_{e^-}]
    \prod_{\substack{f\in E\\ f\ne e}} ([u_{f^+}]+ [u_{f^-}]).
    \]
In \cite{LSZ}, the weights $u_{\gamma}$ and the variable $q$ are restricted to be positive numbers---the case $q=1$ has to be stated separately. 
We observe that the result can be formalized using variables for weights and the parameter $q$, enabling specializations.
For this, we transfer the dependence on $q$ to the operation of sum by introducing  the operation  {\em $q$-sum}, that we denote by $\qcircle$ and define as
\[
a \qcircle b = a + b + (q-1)ab.
\]
 It has the properties that the ordinary sum is recovered with $q=1$, and 
\[
[a+b] = [a] \qcircle [b].
\]
Then, changing $[u_\gamma]$ for $\beta_{\gamma}$, we get
$
d_q(i,j) 
= 
\beta_{i_0 i_{1}} \qcircle \beta_{i_{1}i_{2}} \qcircle\cdots\qcircle \beta_{i_{d-1}i_d} 
$. Thus, the determinant of the matrix of the $d_q(i,j)$ becomes:
\[
    (-1)^{n-1} \sum_{e\in E} \beta_{e^+}\beta_{e^-}
    \prod_{\substack{f\in E\\ f\ne e}} (\beta_{f^+}\qcircle \beta_{f^-}).
    \]

\section{Proof of Theorem \ref{description of the NI in HY} (Existence of a unique non-intersecting path in $\HY$)}\label{technical proofs}

In Theorem  \ref{description of the NI in HY}, we claimed that  for any left-right tree $Y$ whose root has two children, either both right children or  both left children, the network $\HY$ has only one non-intersecting path system, and described it. The proof was postponed. We now proceed to prove this theorem.

Recall that  $V^*(Y)$ is the set of all vertices of $Y$ different from $r$.
 Let us denote with $\Delta$ and $\nabla$ the sets of sources and sinks of $\hemisphere$. 
\begin{align*}
 &\text{Sources: } &&\Delta = 
 \left\{v(i)
 \,|\, i \in \pV(Y) \right\}
 =
 \left\{\M_0(i) M_0(i)
 \,|\, i \in \pV(Y) \right\}
 ,\\
 &\text{Sinks: } &&\nabla =
\left \{e(\arcInput{i}{j})
 \,|\, \arc{i}{j} \in \OriY\right\}
 =
 \left \{M(j,i)M(i,j)
 \,|\, \arc{i}{j} \in \OriY\right\}.
 \end{align*}


 \subsection{Step-disjoint paths and the next-step function associated to a non-intersecting path system.} 

 We introduce two avatars of each  non-intersecting path system $\Lambda=\{\Lambda_1, \ldots, \Lambda_n\}$ of $\HY$.
 
 The first avatar is  the set $\lambda$  of the $n$ paths $\lambda_i$ of $\earring{Y}$ that correspond to the $\Lambda_i$ via the node-arc correspondence.
 Since the $\Lambda_i$ have no node in common, and the nodes of the $\Lambda_i$ are the steps of the $\lambda_i$, the paths $\lambda_i$ are step-disjoint.
The first steps and last steps of the $\lambda_i$ are the elements of the sets $\Delta$ and $\nabla$.
 In particular, the origins and terminus of the $\lambda_i$ are the vertices in the following sets $\Delta_0$ and $\nabla_0$:
  \[
 \Delta_0 = \{\M_0(i)
 \,|\, i \in \pV(Y) \}, \quad
 \nabla_0 = \{M(i,j)
 \,|\, \arc{i}{ j} \in \OriY\}.
 \]

 The second avatar of $\Lambda$ is the \emph{next-step function} $g$ that associates to each step of a path $\lambda_i$, that is not a last step, the next step on the same path.  
 The map $\Lambda \mapsto g$ induces a bijection from  $\NIP(\HY)$ to the set  $\NSF(\earring{Y})$ of all functions fulfilling the following conditions:
 
 \begin{itemize}
 \item[(G1)] The domain $\Domg$ of $g$ is a set of arcs supported on $\earring{Y}$. 
 \item[(G2)] For all $\arc{P}{Q}\in \Domg$, $g(\arcInput{P}{Q})$ is $\arc{Q}{R}$ for some vertex $R \neq  P$.
 \item[(G3)] $g$ is injective.   
 \item[(G4)] $\Delta \subset \Domg$.         
 \item[(G5)] $\nabla \cap \Domg = \emptyset$.
 \item[(G6)] $g(\Domg) \subset \Domg \setminus \Delta \cup \nabla$.
  \end{itemize}

\subsection{Existence of a unique non-intersecting path system in $\hemi{Y}$ }\label{suse:unicity_hemis}

 We show that there is an unique element in $\NSF(\earring{Y})$. 
 
 Let $\Lambda$ be a non-intersecting path system of $\hemisphere$. Let $\lambda$ be the corresponding set of $n$ step-disjoint paths of $\earring{Y}$.
Let $g$ be the next-step function associated to $\lambda$.

For any arc $\arc{P}{Q} \in \arcs{\earring{Y}}$, we define the \emph{flow through $\arc{P}{Q}$} as
\[
\flow(\arcInput{P}{Q})=
\begin{cases}
  1 & \text{ if } \arc{P}{Q} \in \steps{\lambda},\\
  0 & \text{otherwise}.
\end{cases}  
\]
For any vertex $P$ of $\earring{Y}$, we define the \emph{upward  flow from  $P$} as the difference:
\[
\upflow(P) =
\flow(P, \text{parent}(P))
\flow(\text{parent}(P), P).
\]
At any node  $P$ of $\earring{Y}$, we have:
\begin{equation}\label{eq: flow at node}
\upflow(P) - \sum_{Q \text{ child of } P} \upflow(Q) =
\begin{cases}
  +1 & \text{ if } P \in \Delta_0, \ \ \ \ \text{($P$ is an origin)}\\
  -1 & \text{ if } P \in \nabla_0, \ \ \ \ \text{($P$ is a terminus)}\\
  \phantom{+}0 & \text{ otherwise.}
\end{cases}  
\end{equation}

For any vertex $i$ of $Y$ (including the root), let $\Omega_i$ be the set of all arcs $\arc{P}{Q}\in \arcs{\earring{Y}}$ such that $Q$ has contraction $i$. Let $g_i$ be the restriction of $g$ to $\Omega_i$. The sets $\Omega_i$ are the $n+1$ blocks of a partition of $\arcs{\earring{Y}}$. Therefore, $g$ is determined by its restrictions $g_i$. We can work ``locally'' to show the uniqueness of $g$.

We consider now a fixed vertex $i$  of $Y$ different from the root. Let $p$ be its parent,  $j_{-\ell} <_i \cdots <_i j_{-1}$ its left children and  $j_1 <_i \cdots <_i j_m$ its right children. 
We set $\delta=1$ if $\arc{p}{i} \in \OriY$, and $\delta=0$ otherwise. 
To alleviate the notations when working in the local setting ``around $i$'', we introduce an alternative notation for the nodes  incident to arcs of $\Omega_i$. We set: 
\begin{equation}\label{eq: local notation}
\begin{split}
M_0=M_0(i),\quad
\M_0=\M_0(i),\quad
M_\mpp=M(i,p),\quad
\M_\mpp=M(p, i),\\
 M_k = M(i, j_k) \text{ and }  \M_k = M(j_k, i) \text{ for the other values of } k.  
\end{split}
\end{equation}
Note that, with these notations,
$
\delta = 1 \Leftrightarrow \M_\mpp M_\mpp \in \nabla
$.

\begin{figure}
        \includegraphics[width = .49\textwidth, page=35]{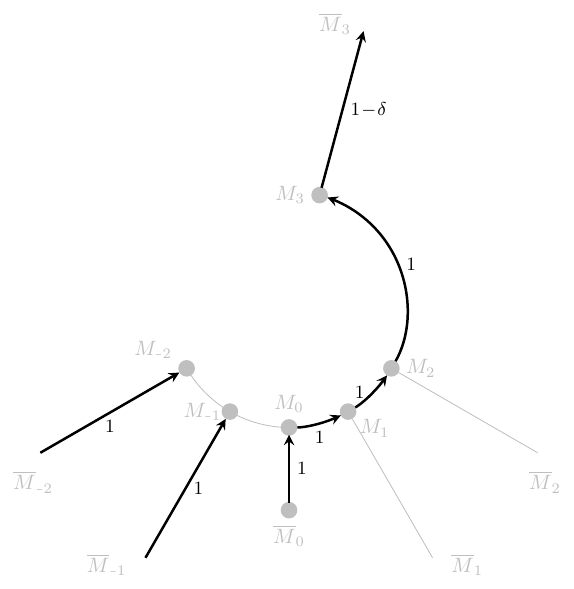}
        \includegraphics[width = .49\textwidth, page=36]{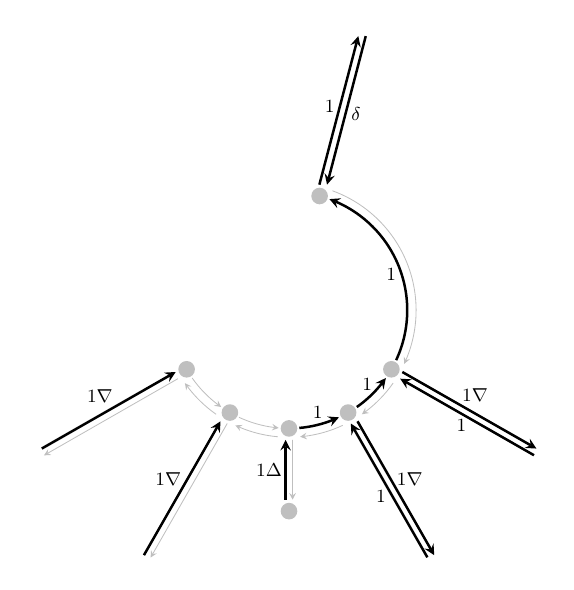}
    \caption{Upward flow $\upflow$ (left) and flow $\flow$ (right) around a vertex $i$ with $2$ left children and 2 right children. 
    Arcs with flow $0$ are left gray, and the flow of all other arcs is indicated by their label. Suffixes indicate 
    first steps and last steps (elements of $\Delta$ and $\nabla$, respectively).}
     \label{fig: upflow}
\end{figure}

The series of  lemmas that follows shows that the upward flow and the flow  are as in Figure \ref{fig: upflow}.  We invite the reader to keep an eye on this Figure when reading the proofs.

\begin{lem}\label{lem: flow from subtree}
  The upward flow satisfies
  \[
  \upflow(\M_k) =
  \begin{cases}
    +1 & \text{ if } k < 0,\\
    0 & \text{ if } 0 < k \le m.
  \end{cases}  
  \]
\end{lem}  

\begin{proof}
Summing \eqref{eq: flow at node} over the set $D$ of all descendants $P$ of $\M_k=\M(i,j_k)$, we get
\[
\upflow(\M_k) = \# D\cap \Delta_0 - \#D\cap \nabla_0.
\]
The contraction maps bijectively   $D\cap \Delta_0$ to the set of vertices in the sub-rooted tree $Y(j_k)$ of $Y$ with root $j_k$.
So $ \# D\cap \Delta_0$ counts the vertices of $Y(j_k)$.

If $k<0$ then $\arc{j_k}{ i} \in \OriY$.
In that case, the set $D\cap \nabla_0$ is in bijection via $M(u, v) \mapsto \{u, v\}$ with the set of edges of $Y(j_k)$.
Then $\#D\cap \nabla_0$ counts the edges of $Y(j_k)$.
Since, as any tree, $Y(j_k)$ has one more vertex than edges,
we get that  $\upflow(\M_k) = +1$.

If $k>0$, then $\arc{j_k}{ i} \not\in \OriY$
and $\# D \cap \nabla_0$  
has one more element than edges in $Y(j_k)$ (the element in excess is $\M_k$), and then $\upflow(\M_k) = 0$.
\end{proof}

\begin{lem}\label{lem: upflow}
  The upward flow satisfies
  \[
  \upflow(M_k)=
  \begin{cases}
    0 & \text{ if } k < 0, \\
    +1 &  \text{ if }0 \le  k \le m,\\
    1-\delta & \text{ if } k=\mpp.
  \end{cases}
  \]
\end{lem}

\begin{proof} 
We apply \eqref{eq: flow at node} at $M_k$ for different values of $k$. We begin at the left children.
\begin{myenum}
    \item If $k = -\ell< 0$, this yields
        $
        \upflow(M_k)  - \upflow(\M_k) = -1 
        $.
        But $\upflow(\M_k)=+1$ by Lemma \ref{lem: flow from subtree}. We deduce  that $\upflow(M_{-\ell}) = 0$.
    \item For $-\ell < k <0$,
        we get 
        $
        \upflow(M_k) - \upflow(M_{k-1}) - \upflow(\M_k) = -1 
        $,
        while $\upflow(\M_k)=+1$ by Lemma \ref{lem: flow from subtree}.
        Then, $\upflow(M_k) = \upflow(M_{k-1})$.

We conclude that $\upflow(M_k)=0$ for all left children, $k < 0$. We consider right children next.

    \item For $k=0$,  we get either 
    \begin{align*}
        \upflow(M_0) -  \upflow(\M_0) &= 0 \quad 
         \text{ if $\ell = 0$, }
        \quad\text{or}\\
        \upflow(M_0) -  \upflow(\M_0) -\upflow(M_{-1}) &= 0 
        \quad \text{ otherwise.}
    \end{align*}
    But in the latter case, $\upflow(M_{-1})=0$ as shown just above. 
    Thus in both cases, we get $\upflow(M_0) =  \upflow(\M_0)$.  
    Since $\upflow(\M_0) = +1$ trivially, we deduce that  $\upflow(M_0)=1$. 
    \item For $0 < k \le m$, we get 
    $\upflow(M_k) - \upflow(M_{k-1}) - \upflow(\M_k) = 0 $, 
    while $\upflow(\M_k)=0$ by Lemma \ref{lem: flow from subtree}.
    Therefore, $\upflow(M_k) = \upflow(M_{k-1})$. 
\end{myenum}

We conclude  that  $\upflow(M_k)=1$ for all right children, $0\le k \le m$.

Finally, in the case $k=\mpp$, we have
$
\upflow(M_{\mpp}) - \upflow(M_m) = -\delta
$ while $\upflow(M_m)=1$ by what precedes.
\end{proof}

\begin{lem}\label{lem: flow vertical}\strut
The flow satisfies
\[
\flow(\M_kM_k) = 
\begin{cases}
\delta & \text{if } k=\mpp \\
1 & \text{otherwise, }
\end{cases}
\;\text{ and }\;
\flow(M_k\M_k) = 
\begin{cases}
0 & \text{if } k\le 0 \\
1 & \text{if } k>0.
\end{cases}
\]
\end{lem}
\begin{proof}
For all $k \le m$ we have
\[\flow(\arcInput{\M_k}{ M_k}) - \flow(\arcInput{M_k }{\M_k}) = \upflow(\M_k).\]
This equation together with the previous lemmas allows us to conclude:
\begin{myenum}
    \item If $k< 0$,  $\upflow(\M_k) = +1$ after Lemma \ref{lem: flow from subtree}, whence $\flow(\arcInput{\M_k}{ M_k})=1$ and  $\flow(\arcInput{M_k }{\M_k})=0$. 
    \item If $0 < k \le m$, $\upflow(\M_k)=0$ whence $\flow(\arcInput{\M_k}{ M_k}) = \flow(\arcInput{M_k }{\M_k})$. Now  $ \flow(\arcInput{M_k }{\M_k})=+1$ since $\arc{M_k}{\M_k} \in \nabla$.
    \item The case $k=0$ is trivial. 
\end{myenum}
Finally, consider $k = \mpp$. We have
\[
\flow(\arcInput{M_{\mpp}}{\M_{\mpp}}) - \flow(\arcInput{\M_{\mpp}}{M_{\mpp}}) = \upflow(M_{\mpp}),
\]
which is $1-\delta$ after Lemma \ref{lem: upflow}, while $\flow(\arcInput{\M_{\mpp}}{ M_{\mpp}})\ge \delta$ by definition of $\delta$. This yields  $\flow(\arcInput{M_{\mpp} }{\M_{\mpp}}) \ge 1$, and thus $\flow(\arcInput{M_{\mpp} }{\M_{\mpp}}) = 1$.
\end{proof}

\begin{lem}\label{lem: flow horizontal}
For all $0 \le k \le m$, the flow satisfies
\[
\flow(\arcInput{M_{k+1}}{M_k})=0
\quad\text{ and }\quad
\flow(\arcInput{M_k}{M_{k+1}})
=\begin{cases}
0 &\text{if } k< 0, \\
1 &\text{if }k \ge 0.
\end{cases}
\]
\end{lem}  

\begin{proof}
By definition of $\upflow$, for all $0\le k\le m$ we have 
\[
\flow(\arcInput{M_k}{M_{k+1}})-\flow(\arcInput{M_{k+1}}{M_k})=\upflow(M_k).
\]
This allows us to conclude for $k\ge 0$, since $\upflow(M_k)=+1$ after Lemma  \ref{lem: upflow}.

For $k <0$, recall $\upflow(M_k)=0$ and thus  $\flow(\arcInput{M_k}{M_{k+1}})=\flow(\arcInput{M_{k+1}}{M_k})$. Our claim is that these are both $0$.
  Let us assume, aiming at a contradiction, that there exists $k < 0$ with $\flow(\arcInput{M_{k+1}}{M_k})=+1$. 
  Consider the smallest such $k$.
  Then $\arc{M_{k+1}}{M_{k}}$ is a step of a path in $\lambda$, but not a last step, while neither $\arc{M_k}{\M_k}$ (that has flow $0$ by Lemma \ref{lem: flow vertical}) nor $\arc{M_k}{M_{k-1}}$ can be its next step (by definition of $k$).
  This gives the desired contradiction. 
\end{proof}  

\begin{lem}\label{lem: graph of gi}
   We have
  \[
  \Domg \cap \Omega_i
  =
  \left\lbrace \arc{\M_k}{M_k} \;|\; 0 \le k \le m \right\rbrace
  \cup \left\lbrace \arc{M_k}{M_{k+1}} \;|\; 0 \le k \le m \right\rbrace.
  \]
  Moreover, for all $0\le k \le m$, we have
  $g(\arcInput{\M_k}{M_k}) = \arc{M_k}{M_{k+1}}$ and $g(\arcInput{M_k}{M_{k+1}}) = \arc{M_{k+1}}{\M_{k+1}}$.
  \end{lem}

\begin{proof}
The arcs in $\Omega_i$ are: the $\arc{\M_k}{M_k}$ for $-\ell \le k \le \mpp$; the $\arc{M_k}{M_{k+1}}$ and $\arc{M_{k+1}}{M_k}$ for $-\ell \le k \le m$. Among them, those that belong to $\Domg$ are those that carry a flow 1 but do not belong to $\nabla$.

The arc $\arc{\M_0}{M_0}$ has flow $1$ since it is in $\Delta$. Trivially $\arc{M_0}{\M_0}$ is not in $\Domg$ since it has no successor in $\hemisphere$.

The arc $\arc{\M_{\mpp}}{M_{\mpp}}$ is not in $\Domg$ since it has flow $\delta$, i.e.~its flow is $1$ if and only it belongs to $\nabla$. 

 All remaining arcs of $\Omega_i$ with flow $1$ are described by Lemmas \ref{lem: flow horizontal} and \ref{lem: flow vertical}. Among them, those in $\nabla$ are the $\arc{M_k }{\M_k}$ for $k > 0$.

 Each $\arc{\M_k}{ M_k}$ for $0 \le k \le m$ has only one successor with flow $1$, which is $M_k M_{k+1}$, whence $g(\arcInput{\M_k}{M_{k}})=\arc{M_{k}}{M_{k+1}}$.  
 Each $\arc{M_{k}}{M_{k+1}}$ admits as successors with flow $1$:   $\arc{M_{k+1}}{\M_{k+1}}$ and (only if $k < m$) $\arc{M_{k+1}}{M_{k+2}}$. 
 But $\arc{M_{k+1}}{M_{k+2}}$ is already the image under $g$ of $\arc{\M_{k+1}}{M_{k+1}}$. Since $g$ is injective, we deduce $g(\arcInput{M_k}{M_{k+1}})=\arc{M_{k+1}}{ \M_{k+1}}$.
\end{proof}

In the previous lemmas, we have completely determined $g$ locally around each vertex $i$ of $Y$ different from the root. We need one last lemma to determine $g$ at the root before closing this section by establishing the uniqueness of non-intersecting path systems of $\hemi{Y}$.

\begin{lem}\label{lem: graph of gr}
  If $r$ has two left children then $\Domg \cap \Omega_r = \emptyset$. 
  If $r$ has two right children, then
$ \Domg \cap \Omega_r = \{\arc{M(a,r)}{ r}, \arc{M(b,r)}{ r}\},$
$ g(\arcInput{M(a,r)}{ r}) = \arc{r}{ M(b,r)}$ and 
$  g(\arcInput{M(b,r)}{r}) = \arc{r}{ M(a,r)}.$
\end{lem}

\begin{proof}
Note that $\Omega_r=\{\arc{M(a,r)}{ r}, \arc{M(b,r)}{ r}\}$. 

If the children of $r$ are left children then $\arc{M(a,r)}{ r}$ and $\arc{M(b,r)}{ r}$ are both in $\nabla$ and thus are not in $\Domg$. On the other hand, if the children of $r$ are right children then $\arc{r}{ M(a,r)}$ is in $\nabla$ but has as only predecessor $\arc{M(b,r)}{ r}$, so $\arc{r}{ M(a,r)}=g(\arcInput{M(b, r)}{ r})$; similarly, mutatis mutandi, $\arc{r}{ M(b,r)}=g(\arcInput{M(a,r)}{r})$.  
\end{proof}

\begin{lem}\label{lem: unique in HY}
  There exists a unique non-intersecting path system of $\HY$.
\end{lem}

\begin{proof}
  Since $\Lambda \mapsto g$ is a bijection from $\NIP(\HY)$ to $\NSF(\earring{Y})$, it is enough to show that $\NSF(\earring{Y})$ has exactly one element.  For any such element $g$,  its restrictions $g_i$ to the subsets $\Omega_i$ are totally determined by Lemmas \ref{lem: graph of gi} and \ref{lem: graph of gr}. 
  This shows that $\NSF(\earring{Y})$ has at most one element.

  To conclude, it remains to 
  check that the function $g$ whose restrictions are the 
  $g_i$ defined in Lemmas \ref{lem: graph of gi} and \ref{lem: graph of gr} belongs indeed to  $\NSF(\earring{Y})$.
  This is done by checking that $g$ fulfills  the
  conditions (G1) to (G6), which is routine. Let us just mention that the injectivity of $g$ (condition (G3)) can be checked locally in each $\Omega_i$. Indeed, for any arc $\arc{P}{Q}$ supported on $\earring{Y}$, all its possible antecedents under $g$ are in the same $\Omega_i$ (namely with $i$ equal to the contraction of $Q$).
\end{proof}

Let $\LU$ be the unique element of $\NIP(\HY)$,
and let $\lu$  the corresponding set of step-disjoint paths of $\earring{Y}$.
For each $i\in V^*(Y)$, let $\LU_i$ be the path of $\LU$ with origin $v(i)$, and let $\lu_i$ be the path of $\lu$ with first step $v(i)$.

To prove Theorem \ref{description of the NI in HY}, we will use the  depth-first search walk of $\earring{Y}$, and the following lemma relating it to the  depth-first search walk of $Y$.

\begin{lem}\label{lem: contraction of DFS}
The contraction of the DFS walk on the earring tree $\earring{Y}$ is the DFS walk on the left-right tree $Y$.
\end{lem}

To prove Theorem \ref{description of the NI in HY},  
we start with  the cyclic DFS walk $w$ of $\earring{Y}$.
The vertices of degree $1$ of $\earring{Y}$ are the $\M(i)$.  
Each $\M_0(i)$  appears exactly once in  $w$, and thus $w$ takes the form
\[
\cdots \M_0(i_1) \cdots \M_0(i_2) \cdots \ldots  \cdots \M_0(i_n) \cdots
\]
where $i_1, i_2, \ldots, i_n$ are the vertices of $Y$ different from the root.
Let $\sigma$ be the cyclic permutation $(i_1i_2 \cdots i_n)$.
The only recoils of the cyclic DFS walk of a tree are at the vertices of degree $1$. 
Thus the only recoils of $w$ are at the $\M_0(i)$, and these recoils are of the form $M_0(i) \M_0(i) M_0(i)$.
As a consequence, the $n$ subwalks 
$
\lambda^\circ_{i} = \M_0(i)  \cdots\M_0(\sigma(i))
$
of the DFS walk $w$ are paths, and each arc supported on $\earring{Y}$ is a step of exactly one of them.

Note that $\lambda^\circ_i$ has first step $\M_0(i) M_0(i)$ and last step $M_0(\sigma(i)) \M_0(\sigma(i))$, since $M_0(i)$ and $M_0(\sigma(i))$ are the only neighbors of the origin and terminus of $\lambda^\circ_i$. 
And, since $\lambda_i^\circ$ has no recoil, it is a path, and can be simply described as the path from $\M_0(i)$ to $\M_0(\sigma(i))$.

Let us call \emph{prefix} of a sequence $x_1 x_2 \cdots x_m$ any sequence $x_1 x_2 \cdots x_k$ for $k \le m$.

\begin{lem}\label{lem: prefix}
  For each $i$, the path $\lu_i$ is a prefix of $\lambda^\circ_{i}$.
\end{lem}

\begin{proof}
  Since the paths in $\lu$ and those in $\lambda^\circ$ have the same first steps,  
it is enough to prove that the ``next step function'' $g$ of $\lu$ is the restriction to $\Domg$ of the next step function $g^\circ$ of $\lambda^\circ$. It is not difficult to describe explicitly the restriction of $g^\circ$ to each $\Omega_i\cap\Domg$.
Consider the case of a vertex $i \neq r$, with $\ell$ left children and $m$ right children. We use again the ``local'' notation $M_k$ and $\M_k$ of \eqref{eq: local notation} for the vertices with contraction $i$ and their neighbors.  Then, for all $k \le m$:
\[
g^\circ(\M_k M_k) = M_k M_{k+1}, \quad\text{ and }\quad 
g^\circ(M_k M_{k+1}) = M_{k+1} \M_{k+1}.
\]
Finally, for the  root,
$
g^\circ(M(a, r) r) = r M(b, r),$ and $g^\circ(M(b, r) r) = r M(a, r). 
$
We conclude  by comparing with Lemmas \ref{lem: graph of gi} and \ref{lem: graph of gr}.
\end{proof}

\begin{lem}\label{lem: middle step}
For each $i \in \pV(Y)$, there exists a unique $h_Y(i) \in \OriY$ such that $e(h_Y(i))$ is a step of $\lambda_i^\circ$. 
Moreover, $\lu_i$ is the path of $\earring{Y}$ whose first step is $v(i)$ and whose last step is $e(h_Y(i))$. 
\end{lem}
\begin{proof}
Recall that $\nabla = \{e(\gamma) \ : \ \gamma\in\OriY\}$.
Since $\lu_i$ is a prefix of $\lambda^\circ_i$, 
the last step of $\lu_i$ is a step of $\lambda^\circ_{i}$.
Therefore,
each path $\lambda^\circ_{i}$ has at least one step in $\nabla$. 
Since the $n$ paths $\lambda^\circ_i$
are step-disjoint and the cardinality of $\nabla$ is also $n$, each of the paths $\lambda^\circ_i$ admits exactly  one step in $\nabla$.
\end{proof}

We can now complete the proof of Theorem \ref{description of the NI in HY}.

\begin{proof}[Proof of Theorem \ref{description of the NI in HY}]
Since $\M_0(i_1), \M_0(i_2),\ldots, \M_0(i_n)$ are the vertices of the form $\M_0(i)$, in their order of visit in the DFS walk of $\earring{Y}$, the sequence $i_1, i_2, \ldots, i_n$ lists the elements of $\pV(Y)$ \emph{in inorder}. 
Indeed,
all vertices with contraction some left child of $i$ are visited before $\M_0(i)$, and all vertices with contraction some right child of $i$ are visited after $\M_0(i)$.   Therefore, $\sigma$ is  the permutation $\sigma_Y$ defined in Theorem \ref{description of the NI in HY}.

For each $i\in\pV(Y)$, the path $\lambda^\circ_{i}$  goes  from $\M_0(i)$ to $\M_0(\sigma(i))$. 
Thus, its contraction is the path of $Y$ from $i$ to $\sigma(i)$: that is path $\mu_i=\pathY{i}{\sigma_Y(i)}$ of  Theorem \ref{description of the NI in HY}.  By Lemma \ref{lem: contraction of DFS},  the cyclic DFS walk of $Y$ decomposes into the paths $\mu_i$. This implies statement \eqref{HY item1}.

The contraction of any walk in $\earring{Y}$ kills all steps except those of the form $e(\arcInput{p}{q})$, that contract into steps $\arc{p}{q}$.
Statement \eqref{HY item2} then follows from Lemma \ref{lem: middle step}.

By Lemma \ref{lem: prefix}, $\lu_i$ is a prefix of $\lambda_i^\circ$. In particular its first step is $v(i)$. Its last step is in $\nabla$, but  $\lambda_i^\circ$ has a  unique step in $\nabla$, namely $e(h_Y(i))$,  by Lemma \ref{lem: middle step}. Therefore, the last step of $\lu_i$ is $e(h_Y(i))$. Note that $h_Y(i)=\arc{j_k}{ j_{k+1}}$. We have obtained that  $\LU_i$ is the path from $v(i)$ to $e(j_k j_{k+1})$, which is, since $i=j_0$, the lifting to $\HY$ of $j_0\cdots j_k j_{k+1}$. This proves  \eqref{HY item3}.

Finally, \eqref{HY item4} can be restated as follows: ``for any arc $\arc{p}{q}$ of $Y$ not in $\OriY$, $e(\arcInput{p}{q})$ is a step of $\lu$ if and only of $q$ is the parent of $p$.'' 
This is checked locally in each  $\Omega_j$ by means of Lemmas \ref{lem: graph of gi} and  \ref{lem: graph of gr}.
\end{proof}

\section{Closing remarks}

While the framework of Choudhury--Khare has proven to be the correct algebraic setting to study distance matrices, our work sets the grounds for a more combinatorial approach to the area. 
We are convinced that catalysts are the right combinatorial objects for studying distance matrices, while route networks provide the natural framework for enumerating catalysts.

To support our idea we mention that a recently released formula for the principal minors of the distance matrix of a tree \cite{Richman} has been elucidated in \cite{Minors} through slight generalizations of the combinatorial objects presented here. Future research directions might explore (i) parametric deformations of this formula and (ii) formulas for the non-principal minors of the matrix.

Another line of research is the study of multiplicative generalizations of the distance matrix of a tree, in the sense of \cite[Thm.~A, case $x\ne0$]{CK19}, \cite[Corollary 2.2 and Thm. 3.3]{YanYeh:2007} and \cite{BS2013}. 
As it stands, these are not explained by Theorem \ref{main: formula indep}.
On a broader level, Choudhury and Khare recently extended in \cite{CK23} their formulas to arbitrary graphs. 
These remain out of the reach of our combinatorial framework.

The determinant of the distance matrix of a graph exemplifies a weighted enumeration of derangements, as evident in Section \ref{sec: catalysts}. Other problems in enumeration of derangements are tackled with our models in \cite{Derangements}.
 
\section*{Acknowledgements}
  
  We dedicate  this  work  to  Ira  Gessel  as  a  token  of  our  gratitude  for  his  lecture  on  the Lindström--Gessel--Viennot  lemma,  which  he  virtually  delivered  at  the  University  of Seville in march of 2022 \cite{GesselSevilla}. The senior authors would also like to express their gratitude to Xavier Viennot for the entrancing memories of his lectures.\medskip
  
  This work has been partially supported by the Grant PID2020-117843GB-I00 funded by MICIU/AEI/ 10.13039/501100011033 and by P2001056 (Junta de Andalucía, FEDER, PAIDI2020).  ÁG was partially funded by the University of Bristol Research Training Support Grant.

\newcommand{\etalchar}[1]{$^{#1}$}

\end{document}